\documentclass[11pt]{article}

\usepackage{amsmath}
\usepackage{amssymb}
\usepackage{lmodern}

\usepackage{bm}
\usepackage{graphicx}
\usepackage{hyperref}
\usepackage{algorithm}
\usepackage{algpseudocode}
\usepackage{tikz}
\usepackage{natbib}
\usepackage{setspace}
\usepackage{caption}
\usepackage{fancyhdr}
\usepackage[margin=1in]{geometry}
\usepackage{titlesec}

\bibpunct[, ]{(}{)}{,}{a}{}{,}%
\usepackage{booktabs}
\usepackage{threeparttable}
\usepackage{etoolbox}
\usepackage[section]{placeins}
\AtBeginEnvironment{algorithmic}{\setstretch{1.2}}


\title{A Pricing-for-Integrality Primal Heuristic for Branch-and-Price: An Application to the Electric Vehicle Routing Problem with Time Windows and Charging Time Slots}

\author{Lukas Eveborn and Elina R\"onnberg \\ 
\small Department of Mathematics, Link\"oping University, 581 83 Link\"oping, Sweden \\
\texttt{lukas.eveborn@liu.se}, \texttt{elina.ronnberg@liu.se}}

\date{}

\titleformat{\subsection}
  {\normalsize\bfseries}{\thesubsection}{1em}{}

\begin{document}

\maketitle

\begin{abstract}
Branch-price-and-cut is the state-of-the-art exact method for solving many types of vehicle routing problems, and it is particularly effective for vehicle routing problems with time windows.
A well-known challenge in branch-price-and-cut is that the generation of columns is guided by information from the linear relaxation of the master problem, with no guarantee that they will be useful from an integer perspective.
As a consequence, high-quality primal solutions are often found only after significant cutting and branching, or the use of primal heuristics.
In this work, based on the ideas of pricing-for-integrality, we propose a primal heuristic for branch-and-price and evaluate it on the electric vehicle routing problem with time windows extended with charging time slots.
The heuristic first uses a destroy-and-repair phase to enrich the column pool and then applies a restricted master heuristic to obtain a primal solution.
In the destroy-and-repair phase, route generation is guided by scaled dual values and heuristic feedback based on previously generated routes, with the aim of producing routes likely to contribute to high-quality integer solutions.
Computational tests on the electric vehicle routing problem with time windows extended with charging time slots show that the proposed heuristic reduces the average optimality gap in the root node by around 60\% compared to a restricted master heuristic and by around 25\%--30\% compared to the implemented benchmark diving heuristic.
\par\noindent\textbf{Funding:} This research was supported by the Swedish Energy Agency within the program FFI, Fordonsstrategisk Forskning och Innovation, under the grant Condore (P2022-00952).
This work was partially supported by the Wallenberg AI, Autonomous Systems and Software Program (WASP) funded by the Knut and Alice Wallenberg Foundation.
\end{abstract}

\section{Introduction}\label{sec:Intro}
Branch-Price-and-Cut (BPC) has, over the past decades, established itself as the state-of-the-art exact method for solving vehicle routing problems, and in particular vehicle routing problems with time windows \citep{ARCHETTI2025, pessoa2020generic}.
A key strength of the BPC algorithm lies in its decomposition of the problem into a master problem and a pricing subproblem.
The pricing problem can often be formulated as a resource-constrained shortest path problem and solved efficiently using a labeling algorithm, allowing column generation to handle the large number of feasible routes without enumerating them explicitly.
The resulting LP relaxations typically provide strong dual bounds, which can be strengthened further through cutting planes.

While column generation is effective for obtaining dual bounds, finding high-quality primal solutions remains challenging.
Primal heuristics are therefore commonly used to search for feasible solutions early in the BPC process.
One common approach is the restricted master heuristic, which solves an integer problem over the columns generated so far.
However, standard pricing generates routes according to their reduced costs with respect to the LP relaxation, and not by their potential contribution to high-quality integer solutions.
Consequently, routes that would be useful for constructing high-quality integer solutions may be missing from the column pool, limiting the performance of the restricted master heuristic.

We therefore propose a pricing-for-integrality primal heuristic that aims to enrich the column pool with routes that are more likely to contribute to high-quality integer solutions before applying a restricted master heuristic over the expanded column pool.
The heuristic begins by constructing several initial integer solutions from the existing column pool, allowing different starting points to be explored without requiring an incumbent solution.
Each initial solution is partially destroyed using a spatially aware procedure, after which repair iterations are performed to generate new routes.
Central to these iterations is a repair-pricing scheme inspired by the pricing-for-integrality heuristic IPColGen \citep{maher2023integer}.

IPColGen was developed to generate columns that are likely to contribute to high-quality feasible solutions.
In IPColGen, a partial integer solution is first constructed by removing columns from the incumbent solution.
The partial solution is then repaired iteratively through modified pricing.
The pricing objective combines the original cost, a scaled dual contribution, and heuristic feedback derived from the partial solution and previously generated repair columns.
In IPColGen, repair pricing places particular emphasis on generating columns that complement the partial solution.
After the repair-pricing iterations, a repair sub-MIP is solved over the generated columns, with the partial solution fixed, to search for an improved feasible integer solution.
Our heuristic retains the three components of the modified pricing objective but uses them with a broader focus on column-pool enrichment.
To do so efficiently, we exploit the ability of a labeling algorithm to solve the repair pricing problem heuristically and return several routes in each iteration.
The heuristic feedback is also adapted to guide the set of generated routes more gradually by responding less strongly to deviations during the repair process.

Together, these components allow the heuristic to explore different parts of the solution space and efficiently build a broad pool of promising routes from an integer perspective, giving the restricted master heuristic more opportunities to find a high-quality primal solution.

To evaluate the heuristic, we introduce the Electric Vehicle Routing Problem with Time Windows and Charging Time Slots (EVRPTW-CTS), which extends the EVRPTW \citep{Schneider2014} with capacitated charging resources available during specific time slots.
The problem is motivated by practical considerations regarding the need to coordinate shared charging infrastructure for heavy-duty electric vehicles, both to avoid waiting times and to reduce power peaks.
In a joint project with the truck manufacturer Scania, we study bookable charging time slots as a means of addressing this challenge, and simulate other vehicles’ use of the chargers by blocking selected time slots in our instances.
This contrasts with most of the EVRP literature, which assumes that charging resources are abundant and always available, a simplifying assumption that is often unrealistic in practice.
Some previous studies consider limited charging resources \citep{LAM2022105870,froger2022electric}.
However, an important difference in our instances is that the charging limitations have a substantial impact on the resulting routes.

The main contributions of this work are:
\begin{itemize}
    \item We propose a pricing-for-integrality primal heuristic for branch-and-price, tailored to vehicle routing problems.
    It enriches the column pool with routes likely to contribute to high-quality integer solutions before applying a restricted master heuristic to the expanded pool.
    \item We introduce the EVRPTW-CTS, which extends the EVRPTW with capacitated charging resources available for booking during specified time slots, and provide benchmark instances in which these restrictions have a clear effect on the solutions.
    \item We computationally evaluate the proposed heuristic against restricted master and diving heuristics on the EVRPTW-CTS instances and analyze its use at different stages of the BPC process.
    \item The code and instances are publicly available and can be found at \url{https://gitlab.liu.se/lukev81/bpc-for-evrptw-cts-scip}.
\end{itemize}

The remainder of the paper is structured as follows.
In Section \ref{sec:Background}, relevant background and related work regarding BPC, pricing for integrality, and electric vehicle routing problems is presented.
This is followed by a problem statement for the EVRPTW-CTS in Section \ref{sec:ProblemFormulation}.
In Section \ref{sec:HeuristicPricing}, we provide a detailed description of the heuristic, including algorithmic details and design choices.
The BPC-framework used to solve the EVRPTW-CTS and analyze the performance of the proposed heuristic is described in Section \ref{sec:BPC}.
Subsequently, computational results and insights are presented in Section \ref{sec:Results}.
Finally, conclusions and future work are discussed in Section \ref{sec:Conclusion}.

\section{Background and Related Work}\label{sec:Background}
This section provides an overview of the relevant literature related to our work, focusing on branch-price-and-cut, pricing-for-integrality and electric vehicle routing problems.

\subsection{Branch-Price-and-Cut}
Branch-price-and-cut combines column generation, branching, and cutting planes to solve an extended formulation with too many variables to enumerate explicitly.
At each node of the search tree, the linear relaxation is solved over a restricted set of columns.
The resulting dual values are used in a pricing problem to identify columns with negative reduced cost, which are added to the Restricted Master Problem (RMP).
When no such columns remain, the linear relaxation has been solved to optimality; branching is then applied if the solution is fractional.
Cutting planes can be added to strengthen the relaxation.

In vehicle routing, the pricing problem is commonly formulated as an Elementary Shortest Path Problem with Resource Constraints (ESPPRC).
It is typically solved by a labeling algorithm that tracks resource consumption and uses dominance rules to discard labels that cannot lead to a better route.
For further details on BPC and labeling algorithms, we refer to \cite{G-2024-36,uchoa2024optimizing} and \cite{irnich2005shortest}, respectively.
However, a few words about primal heuristics in BPC are warranted, as they are central to our work.

When it comes to primal heuristics in BPC, it has turned out to be more complex than for branch-and-bound.
\cite{G-2024-36} groups primal heuristics for BPC into five categories: restricted master heuristics, complementary pricing, rounding and sub-MIP heuristics, diving heuristics, and large-neighbourhood search.
Restricted master heuristics and complementary pricing are typically easier to implement, whereas diving heuristics, rounding and sub-MIP heuristics, and large-neighbourhood search are generally more sophisticated and may obtain higher-quality solutions by using information from deeper in the search tree or through reoptimization.

When it comes to vehicle routing problems and BPC, the most common primal heuristics are based on restricted master heuristics and diving heuristics (see e.g. \cite{joncour2010column}, \cite{spliet2015discrete} and \cite{fink2019column}).
A restricted master heuristic and diving heuristics are also what is implemented in VRPSolver \citep{pessoa2020generic}, a state-of-the-art branch-price-and-cut framework for solving a broad range of vehicle routing problems by exact methods.
There the authors use the diving heuristics of \cite{sadykov2019primal}.
Overall, however, branch-and-price-based primal heuristics remain relatively uncommon in the vehicle routing literature.

\subsection{Pricing-for-Integrality}\label{subSec:PricingForIntegralityBackground}
A recent primal heuristic for BPC is IPColGen, proposed by \cite{maher2023integer}.
IPColGen is a column-generation-based matheuristic embedded in a large-neighbourhood search framework with destroy and repair phases.
Its main contribution is a repair pricing scheme designed to generate columns that are likely to contribute to high-quality integer solutions of set covering, packing, and partitioning problems.
The scheme builds on the global optimality conditions of \cite{larsson2006global} and the integer programming column-generation ideas of \cite{zhao2020integer}, but generates a collection of repair columns over successive repair iterations rather than greedily replacing one column at a time.

IPColGen starts from the current incumbent and removes some of its columns to form a partial solution.
Repair pricing then generates columns intended to complement the remaining columns.
The modified pricing objective combines the original cost, an LP dual contribution scaled by $\gamma$, and row-specific penalties based on the partial solution and previously generated repair columns.
A large static penalty discourages repair columns from using packing rows already satisfied at equality by the partial solution.
Dynamic penalties encourage coverage of uncovered covering rows and discourage excessive use of packing rows not filled by the partial solution.
The dynamic penalties are adjusted multiplicatively according to the contributions of the columns generated so far.
Each repair iteration generates one column, which is added to the RMP and a repair sub-MIP.
After the repair iterations, the sub-MIP searches for an improved solution with the partial solution fixed.

A brief motivation for the modified pricing objective in IPColGen can be obtained from a Lagrangian-dual perspective; for a detailed theoretical discussion, we refer to the IPColGen paper \citep{maher2023integer}.
The original costs determine the quality of an integer solution, whereas reduced costs based on the LP dual values guide column generation towards solving the LP relaxation.
Because of the duality gap, the columns with the most attractive reduced costs do not necessarily form a high-quality integer solution, and identifying an optimal integer solution based on reduced costs may require the enumeration of a large number of columns; see, e.g., \cite{ronnberg2019integer}.
As an alternative, \cite{zhao2020integer} used the results of \cite{larsson2006global} to propose scaling of the contributions of the dual values to tilt the hyperplanes of the Lagrangian dual function as a heuristic strategy for generating columns for integer solutions.
The repair pricing scheme of IPColGen is a further development of this idea.

IPColGen is invoked adaptively after column-generation tailing-off has been detected, primarily at the root node but also at selected tree depths.
Computational experiments on both structured decomposable problems and general instances from MIPLIB 2017 show improvements in the final optimality gap and primal integral compared with a restricted master heuristic, particularly for instances with large optimality gaps when IPColGen is first called.

\subsection{Electric Vehicle Routing Problems}\label{subSec:ElectricVehicleRoutingProblems}
Electric Vehicle Routing Problems (EVRPs) extend traditional vehicle routing problems with constraints on battery capacity and charging.
Overviews are provided by \cite{erdelic2019survey} and \cite{KUCUKOGLU2021107650}, and a detailed background by \cite{montoya2016electric}.
The Electric Vehicle Routing Problem with Time Windows (EVRPTW), introduced by \cite{Schneider2014}, requires capacitated vehicles to serve customers within specified time windows while accounting for their limited driving range and possible visits to charging stations.

The formulation of \cite{Schneider2014}, like many subsequent formulations, assumes unlimited charging-station capacity.
In practice, limited numbers of chargers and available power can cause queues and disrupt route schedules.
Some studies address this by explicitly modelling waiting times at charging stations \citep{keskin2021simulation,keskin2019electric,kullman2021electric}.
Others coordinate routing and charging to avoid conflicts: \cite{ding2015conflict} adjust the schedules of routes generated without charger-capacity constraints, \cite{LAM2022105870} combine Dantzig--Wolfe and logic-based Benders decomposition to solve instances with up to 100 customers to optimality, and \cite{froger2022electric} use a matheuristic that accounts for charger capacity and addresses instances with up to 320 customers.

Under the experimental conditions in these studies, limited charging capacity has generally had a modest effect on solutions.
\cite{froger2022electric} report an average objective-value increase of 0.29\% relative to unlimited capacity, \cite{LAM2022105870} find that capacity constraints affect few instances, and \cite{ding2015conflict} observe that charging conflicts become consequential mainly when capacity is nearly insufficient.
However, these findings may partly reflect assumptions that do not fully capture real-world operations.
Many chargers are public and shared among different users, so vehicles belonging to one fleet cannot freely coordinate their charging activities.
Moreover, commercial vehicles typically operate under tight schedules and high utilization rates.
Unplanned waiting at a charger may therefore cause delays that propagate through the schedule, affecting delivery reliability and operational efficiency.

To study charging capacity under these conditions, we consider bookable time slots at charging stations.
This concept, explored in collaboration with the truck manufacturer Scania, allows vehicles to reserve charging times in advance and makes their schedules less dependent on uncertain charger availability.
We use it to investigate how limited and shared charging capacity affects routing decisions when competition for charging resources is more pronounced.
The following section defines the resulting EVRPTW-CTS and presents its mathematical formulation.

\section{Problem Formulation for the EVRPTW-CTS}\label{sec:ProblemFormulation}
In this section the problem is formally defined, and an extended formulation is presented.
A compact formulation of the problem can be found in Appendix \ref{sec:CompactModel}.

\subsection{Problem statement}
In the Electric Vehicle Routing Problem with Time Windows and Charging Time Slots (EVRPTW-CTS), a fleet of homogeneous electric vehicles must serve a set of customers.
The vehicles are based at a central depot and have a limited battery capacity as well as a limited carrying capacity.
Each customer has a demand that must be fulfilled and a time window in which the service must start.
A vehicle is allowed to wait if it arrives before the time window opens, but arriving after the time window closes is not permitted.
To be able to serve all customers, the vehicles can visit charging stations to recharge their batteries.
Motivated by real-world applications, described in Section \ref{subSec:ElectricVehicleRoutingProblems}, charging stations can only be used during specific time slots.
Even if a charger is not used for the entire time slot, it cannot be used by another vehicle during that time slot.
Outside the time slots, the chargers are unavailable, but it is allowed to arrive early and wait until the time slot begins.
The vehicles can partially recharge their batteries, and the recharging is modelled to be linear.
The objective is to minimize the total distance driven by the vehicles while ensuring that all customers are served within their time windows.
This while making sure that the vehicles do not run out of battery or exceed their carrying capacity.

\subsection{Extended formulation of the EVRPTW-CTS}
Consider a directed graph $G = (N, A)$, where $N$ is the set indexing the nodes and $A$ is the set indexing the arcs.
Let $V \subseteq N$ be the set indexing the customers and $F \subseteq N$ be the set indexing the charging time slots, where each charging time slot is defined by its start time, end time, and the location of the charging station.
Let $K$ denote the set of available vehicles and $P$ denote the set of feasible routes.
A route $p \in P$ starts and ends at the depot, visits each customer or charging time slot at most once, satisfies vehicle battery and capacity constraints, and respects customer as well as charging time windows.
For each route $p \in P$, let $c_p$ denote its cost and let $a_{ip}$ be a binary parameter equal to 1 if node $i \in V \cup F$ is visited by route $p$, and 0 otherwise.
The binary variable $\lambda_p$ indicates whether route $p \in P$ is selected or not.
An extended formulation of the problem can then be stated as follows:
  \begin{subequations}\label{eq:extended_model}
  \begin{align}
  \text{minimize} \quad
      & \sum_{p \in P} c_p \lambda_p  & \label{extended_model:objective} \\
  \text{subject to} \quad
      & \sum_{p \in P}  a_{ip} \lambda_p \ge 1, & i \in V, \label{extended_model:customer_visits} \\
      & \sum_{p \in P}  a_{ip} \lambda_p \le 1, & i \in F, \label{extended_model:charging_stations_visit} \\
      & \sum_{p \in P}  \lambda_p \le |K|, &  \label{extended_model:vehicle_capacity} \\
      & \lambda_p \in \{0,1\}, & p \in P.
  \end{align}
  \end{subequations}

The objective function \eqref{extended_model:objective} minimizes the total cost of the selected routes.
Constraints \eqref{extended_model:customer_visits} ensure that each customer is visited at least once, while constraints \eqref{extended_model:charging_stations_visit} ensure that each charging time slot is visited at most once.
In addition, constraint \eqref{extended_model:vehicle_capacity} ensures that the number of used vehicles does not exceed the available fleet size.

The restricted master problem (RMP) is the LP relaxation of the extended formulation with only a subset of the routes in $P$ included.
The associated pricing problem is an ESPPRC whose objective is to find $\min_{p\in P}\bar{c}_p$, where $\bar{c}_p$ is the reduced cost of route $p$.
Let $\sigma_i$ be the dual value of the constraint for node $i \in V \cup F$ in the corresponding RMP, and let $\mu$ be the dual value of the vehicle limit constraint.
The reduced cost of a route $p$ is then $\bar{c}_p=c_p-\sum_{i\in V}\sigma_i a_{ip}-\sum_{i\in F}\sigma_i a_{ip}-\mu$, with $\sigma_i \ge 0$ for $i \in V$, $\sigma_i \le 0$ for $i \in F$, and $\mu \le 0$.

\section{Pricing for Integrality Heuristic}\label{sec:HeuristicPricing}
In this section, we present the proposed pricing-for-integrality heuristic, a primal heuristic for branch-and-price tailored to vehicle routing problems.
The heuristic combines a destroy-and-repair phase that expands the column pool with a restricted master heuristic that searches for a primal solution over the expanded pool.

The heuristic addresses a limitation of standard column generation.
Standard column generation uses dual information from the current RMP solution to generate columns that are attractive from the perspective of the LP relaxation.
However, such columns are not necessarily compatible, or attractive, in an integer solution, which may limit the quality of the solutions obtainable from the generated column pool.

To address this limitation, the proposed heuristic uses a destroy-and-repair phase to enrich the column pool through a modified pricing process.
This phase can be briefly summarized as follows.
In each iteration, an initial integer solution is partially destroyed, after which pricing is performed on the resulting repair graph.
The objective coefficients of the pricing problem combine the original arc costs ands scaled dual values with feedback from columns generated earlier in the current repair process.
Consequently, the generated columns are encouraged to complement the remaining part of the initial integer solution and collectively satisfy the master problem constraints associated with the nodes in the repair graph, while preserving the guidance provided by the dual values.

The remainder of this section first introduces the problem setting and notation, then presents the algorithmic scheme, and finally describes the individual components of the heuristic in detail.

\subsection{Problem Setting and Notation}\label{subsec:ProblemSettingAndNotation}
For the purpose of describing the heuristic, we consider a vehicle routing problem with homogeneous vehicles defined on a directed graph $G=(N,A)$ with one depot node.
We assume that each node $i \in N$ is either associated with a master problem constraint of the form set covering, set partitioning, or set packing, or part of no master problem constraint at all.
For example, in the EVRPTW-CTS, customer nodes are associated with set covering constraints and charging time-slot nodes with set packing constraints.
Each column generated by the pricing problem represents a route in $G$ that visits a subset of the nodes in $N$ and starts and ends at the depot node.
Each arc $(i,j) \in A$ has a cost $c_{ij}$, and the cost of a column is the sum of the costs of the arcs traversed by its corresponding route.
The objective of the master problem is to select a feasible set of columns that minimizes the total cost.

At any point during the branch-and-price process, the current RMP contains a subset of the columns in the complete master problem and may additionally contain valid inequalities and other constraints introduced during the solution process.
Let $C$ denote the current pool of generated columns, which may include columns that are not currently present in the RMP.
The pricing problem uses dual information obtained from a dual solution to the current RMP.
Let $\boldsymbol{\pi}$ denote the corresponding vector of dual values associated with the original master constraints, valid inequalities, and any additional constraints included in the RMP.
For each arc $(i,j) \in A$, $\pi_{ij}$ incorporates all dual contributions from $\boldsymbol{\pi}$ that can be attributed to traversing that arc.

The pricing problem seeks a feasible route with minimum reduced cost with respect to $\boldsymbol{\pi}$.
We assume that the pricing problem can be solved using a labeling algorithm, as is common for vehicle routing problems.

Although presented in the context of vehicle routing, the heuristic is also applicable to problems such as crew scheduling, where pricing can be solved using a state-based algorithm.

\subsection{Algorithmic Overview}
The overall scheme of the heuristic is presented in Algorithm~\ref{alg:heuristic_pricing_for_integrality}.
The heuristic can be applied at any point during the branch-and-price process at which a column pool $C$ and a vector of dual values $\boldsymbol{\pi}$ are available.
We next describe the inputs and the steps of the algorithm.

In addition to \(C\) and \(\pi\), the heuristic takes several parameters as input.
The parameters \(I_{\mathrm{start}}\) and \(I_{\mathrm{destroy}}\) specify the number of initial solutions to construct and the number of destroy-and-repair iterations performed for each initial solution, respectively.
The parameter \(\rho\) denotes the fraction of the solution destroyed in each destroy-and-repair iteration.
Finally, the parameters \(\gamma\) and \(\beta\) balance the influence of the dual values and the heuristic weights in the modified pricing objective.

For each of the \(I_{\mathrm{start}}\) start iterations, an initial integer solution \(S^{\mathrm{start}}\) is constructed from the current column pool, as described in Section~\ref{subsec:CreatingInitialSolutions}.
The heuristic then performs \(I_{\mathrm{destroy}}\) destroy-and-repair iterations starting from this solution.
In each such iteration, a fraction \(\rho\) of the columns in \(S^{\mathrm{start}}\) is removed, leaving a partial solution to be repaired.
Based on the removed columns, a repair graph \(G^{\mathrm{rep}}=(N^{\mathrm{rep}},A^{\mathrm{rep}})\), with \(N^{\mathrm{rep}}\subseteq N\) and \(A^{\mathrm{rep}}\subseteq A\), is constructed as described in Section~\ref{subsec:CreateRepairGraph}.
The construction of the repair graph also determines the number of repair iterations \(I_{\mathrm{repair}}\) based on the number of removed columns.
Hence, \(I_{\mathrm{repair}}\) may vary across destroy-and-repair iterations.

The repair process iteratively solves the repair pricing problem to generate columns that complement the partial solution.
Since a labeling algorithm is used, multiple columns can be generated in each repair iteration.
Moreover, as the overall procedure is heuristic, optimality of the repair pricing solutions is not required.
The repair pricing problem can therefore be solved heuristically, limiting the computational effort per iteration and allowing a relatively large number of columns to be generated within a limited amount of time.

To provide feedback across repair iterations, let $\boldsymbol{U}_r$ denote the vector of cumulative node-visit measures after iteration $r$ of the current repair process.
Its entries summarize how frequently nodes have been visited by the columns generated so far.

The first repair iteration is treated separately because the heuristic weights are based on columns generated in earlier iterations of the current repair process, and no such columns are initially available.
The repair pricing problem is therefore solved using the standard reduced arc costs $\bar{\boldsymbol{c}}$, computed from the original arc costs and the dual values $\boldsymbol{\pi}$.
This produces a set of columns $C^{\mathrm{rep}}_1$.
These columns are added to the column pool \(C\) and used to initialize the vector of cumulative node-visit measures \(\boldsymbol{U}_1\).

In each subsequent repair iteration \(r=2,\ldots,I_{\mathrm{repair}}\), the pricing objective is modified based on the dual values \(\pi\), the cumulative node-visit measures \(\boldsymbol{U}_{r-1}\), the current iteration \(r\), and the total number of repair iterations \(I_{\mathrm{repair}}\).
This modification provides feedback based on the difference between each node's cumulative visit measure and its expected value at the current stage of the repair process.
Depending on the type of the associated master problem constraint, the modified objective encourages or discourages columns that visit the node, with the aim of generating columns that collectively satisfy the constraints of the repair problem, as described in Section~\ref{subsec:ModificationOfPricingObjective}.
Solving the repair pricing problem with the modified objective coefficients \(\bar{\boldsymbol{c}}_r^{\mathrm{heur}}\) produces a new set of columns \(C^{\mathrm{rep}}_r\).
These columns are added to \(C\) and used to update the cumulative node-visit measures \(\boldsymbol{U}_r\).

After all start and destroy-and-repair iterations have been completed, a restricted master heuristic is applied to the expanded column pool to obtain a primal integer solution \(S^{\mathrm{primal}}\).

\begin{algorithm}
\begin{algorithmic}[1]
\caption{Pseudocode of the pricing-for-integrality heuristic}\label{alg:heuristic_pricing_for_integrality}
\Require Column pool $C$ and dual values $\boldsymbol{\pi}$
\Require Number of initial solutions $I_{\mathrm{start}}$ and destroy-and-repair iterations $I_{\mathrm{destroy}}$
\Require Destroy fraction $\rho$ and weighting parameters $\gamma$ and $\beta$
\Ensure Primal integer solution $S^{\mathrm{primal}}$
\Procedure{PricingForIntegralityHeuristic}{$C,\boldsymbol{\pi},I_{\mathrm{start}},I_{\mathrm{destroy}},\rho,\gamma,\beta$}
    \For{$m \gets 1$ \textbf{to} $I_{\mathrm{start}}$}
        \State $S^{\mathrm{start}} \gets \textsc{CreateInitialSolution}(C)$
        \For{$n \gets 1$ \textbf{to} $I_{\mathrm{destroy}}$}
            \State $G^{\text{rep}}, I_{\mathrm{repair}} \gets \textsc{ConstructRepairGraph}(S^{\mathrm{start}}, \rho)$
            \State $C^{\text{rep}}_1 \gets \textsc{SolveRepairPricingProblem}(G^{\text{rep}}, \boldsymbol{\bar{c}})$
            \State $C \gets C \cup C^{\text{rep}}_1$
            \State $\boldsymbol{U_1} \gets \textsc{UpdateCumulativeVisits}(\boldsymbol{0}, C^{\text{rep}}_1)$
            \For{$r \gets 2 $ \textbf{to} $I_{\mathrm{repair}}$}
                \State $\boldsymbol{\bar{c}_{r}^{\text{heur}}} \gets \textsc{ModifyPricingObjective}(G^{\text{rep}}, \boldsymbol{\pi}, \boldsymbol{U_{r-1}}, r, I_{\mathrm{repair}}, \gamma, \beta)$
                \State $C^{\text{rep}}_r \gets \textsc{SolveRepairPricingProblem}(G^{\text{rep}}, \boldsymbol{\bar{c}}_{r}^{\text{heur}})$
                \State $C \gets C \cup C^{\text{rep}}_r$
                \State $\boldsymbol{U_r} \gets \textsc{UpdateCumulativeVisits}(\boldsymbol{U_{r-1}}, C^{\text{rep}}_r)$
            \EndFor
        \EndFor
    \EndFor
    \State $S^{\mathrm{primal}} \gets \textsc{RestrictedMasterHeuristic}(C)$
    \State \Return $S^{\mathrm{primal}}$
\EndProcedure
\end{algorithmic}
\end{algorithm}

\subsection{Creating Initial Solutions}\label{subsec:CreatingInitialSolutions}
Each start iteration of the heuristic begins by constructing an initial integer solution, $S^{\mathrm{start}}$.
Since the subsequent destroy-and-repair procedure is applied separately to each initial solution, diversity among them allows the heuristic to explore different parts of the solution space.

In our implementation, the initial solutions are constructed by repeatedly applying a restricted master heuristic to the column pool \(C\).
To ensure that an initial integer solution can be constructed even when the column pool is insufficient, we add high-cost unit columns for the set covering and set partitioning constraints.
To promote diversity, a small random perturbation is applied to the column costs before each initial solution is constructed.
For each column \(p\in C\), the modified objective coefficient is given by \(c_p^{\mathrm{modified}}=c_p(1+\xi_p)\), where \(\xi_p\sim\mathcal{U}(-\epsilon,\epsilon)\).
The parameter \(\epsilon\) controls the magnitude of the perturbation and is set to a small value.
We used \(\epsilon=0.05\) in our experiments.

Since the initial solutions only serve as starting points for the subsequent destroy-and-repair procedure, limited computational effort is devoted to their construction.
The restricted master heuristic is therefore terminated when either a prescribed time limit or an optimality-gap limit is reached, and the best solution found is used as the initial solution.
Based on preliminary experiments, these limits were set to 10 seconds and \(5\%\), respectively.

\subsection{Constructing the Repair Graph}\label{subsec:CreateRepairGraph}
Given an initial solution $S^{\mathrm{start}}$, the repair graph $G^{\text{rep}}$, is constructed in a two-step process.
First, a partial integer solution is created by removing a subset of the columns from the initial solution.
Let $N^{\text{partial}}$ denote the set of nodes visited by the columns in the partial solution.
The repair graph is then constructed from the nodes that are not in the partial solution, together with the depot node, which is denoted as node 0.
Hence, $N^{\mathrm{rep}} = \left(N \setminus N^{\mathrm{partial}}\right) \cup \{0\}$ and $A^{\mathrm{rep}} = \left\{(i,j)\in A : i,j\in N^{\mathrm{rep}}\right\}$.

Columns are selected for removal using a randomized procedure based on the spatial proximity of the routes they represent.
First, one route is selected uniformly at random from $S^{\mathrm{start}}$ and designated as the seed route $p^{\mathrm{seed}}$.
For each remaining route $p$, let $d(p,p^{\mathrm{seed}})$ denote the Euclidean distance between the centroid of $p$ and that of $p^{\mathrm{seed}}$.
The probability of selecting route $p$ for removal is then given by $\bigl(1/d(p,p^{\mathrm{seed}})\bigr)\big/\sum_{q\in S^{\mathrm{remaining}}}\bigl(1/d(q,p^{\mathrm{seed}})\bigr)$, where $S^{\mathrm{remaining}}$ denotes the set of routes that have not yet been removed.
Thus, routes whose centroids are closer to that of the seed route are more likely to be removed.
Since geographically close nodes are more likely to be combined in a low-cost route, removing nearby routes gives repair pricing an opportunity to recombine their nodes in new ways.

In each iteration of the selection procedure, exactly one route is sampled and removed from $S^{\mathrm{remaining}}$.
The selection probabilities are then recomputed before the next route is sampled.
With $\rho$ as the destroy fraction, the procedure continues until a total of $\left\lfloor\rho|S^{\mathrm{start}}|\right\rfloor$ routes have been removed, including the seed route.
The routes that remain constitute the partial solution.

\subsection{Modification of Pricing Objective and Update of Cumulative Visits}\label{subsec:ModificationOfPricingObjective}
Modifying the pricing objective is a central component of the heuristic.
It guides the search towards columns that are likely to be part of high-quality integer solutions.

In a standard column generation setting, the pricing problem is formulated with the goal of finding columns with negative reduced cost based on the dual values from the RMP.
The reduced cost contribution for an arc $(i,j) \in A^{\text{rep}}$ is typically calculated as $\bar{c}_{ij} = c_{ij} - \pi_{ij}$, where $c_{ij}$ is the cost of the arc and $\pi_{ij}$ incorporates all dual contributions associated with traversing arc $(i,j)$.
Rather than using these standard reduced-cost contributions, the modified pricing objective combines three components: (1) original arc costs, (2) scaled dual contributions, and (3) heuristic weights.

The heuristic weights provide feedback based on the columns generated during the current repair process.
Their purpose is to encourage the generated columns to collectively satisfy the master problem constraints associated with the nodes in the repair pricing problem.
For each node, the weight is based on the difference between its expected cumulative visit measure and its cumulative visit measure in the columns generated so far.

Let $U_{jr}$ denote the cumulative visit measure of node $j \in N^{\mathrm{rep}}$ after repair iteration $r$, with $U_{j0}=0$.
It is updated according to $U_{jr} = U_{j(r-1)} + |\{c \in C_r^{\mathrm{rep}} : c \text{ visits } j\}|/|C_r^{\mathrm{rep}}|$.
Thus, each iteration adds the fraction of its generated columns that visit node \(j\) to the cumulative measure.

Before iteration $r$, the expected cumulative visit measure is $\frac{r-1}{I_{\mathrm{repair}}}$, where $I_{\mathrm{repair}}$ is the total number of repair iterations in the current destroy-and-repair iteration.
It therefore increases by $\frac{1}{I_{\mathrm{repair}}}$ per iteration, from zero before the first iteration to a final value of one after all repair iterations.

The difference $\frac{r-1}{I_{\mathrm{repair}}}-U_{j(r-1)}$ indicates whether the cumulative visit measure of node $j$ is below or above its expected value before iteration $r$.
How this difference should affect the repair pricing objective depends on the type of the associated master problem constraint.
For a set covering constraint, insufficient visits should make columns visiting the node more attractive, whereas excessive visits require no correction.
For a set packing constraint, excessive visits should make columns visiting the node less attractive, whereas insufficient visits require no correction.
For a set partitioning constraint, insufficient visits should make such columns more attractive and excessive visits should make them less attractive.

To scale this feedback to the cost structure of the instance, we use the average arc cost in the full graph, defined as $c^{\mathrm{avg}}=\sum_{(i,j)\in A}c_{ij}/|A|$.
The heuristic weight of node $j \in N^{\mathrm{rep}}$ in repair iteration $r$ is then defined as

\begin{equation}
h_{jr}
=
\begin{cases}
\displaystyle
c^{\mathrm{avg}}
\max\left\{0,\frac{r-1}{I_{\mathrm{repair}}}-U_{j(r-1)}\right\},
& \text{if $j$ is associated with a set covering constraint},\\[1ex]
\displaystyle
c^{\mathrm{avg}}
\min\left\{0,\frac{r-1}{I_{\mathrm{repair}}}-U_{j(r-1)}\right\},
& \text{if $j$ is associated with a set packing constraint},\\[1ex]
\displaystyle
c^{\mathrm{avg}}
\left(\frac{r-1}{I_{\mathrm{repair}}}-U_{j(r-1)}\right),
& \text{if $j$ is associated with a set partitioning constraint},\\[1ex]
0,
& \text{if $j$ is not associated with any such constraint}.
\end{cases}
\label{HPFI:heuristic_weight}
\end{equation}

The resulting heuristic weights constitute the node-specific feedback component of the repair pricing objective.
They are combined with the original arc costs and the scaled dual contributions, where the previously introduced parameters $\gamma$ and $\beta$ control the scaling.
The modified objective coefficient of arc $(i,j)\in A^{\mathrm{rep}}$ in repair iteration $r$ is then defined as
\begin{equation}
        \bar{c}_{ijr}^{\mathrm{heur}} = c_{ij}-\gamma\pi_{ij}-\beta h_{jr}. \label{HPFI:modified_objective}
\end{equation}
Following the motivation of \citet{maher2023integer} discussed in Section~\ref{subSec:PricingForIntegralityBackground}, $\gamma$ scales the dual contribution to moderate its influence relative to the original arc cost.
The parameter $\beta$ independently controls the strength of the cumulative-visit feedback.

\section{Branch-Price-and-Cut Framework}\label{sec:BPC}
In order to solve the Electric Vehicle Routing Problem with Time Windows and Charging Time Slots (EVRPTW-CTS) and evaluate the performance of the proposed pricing-for-integrality heuristic, a branch-price-and-cut algorithm was implemented using SCIP \citep{bolusani2024scip} as the underlying framework.
This section describes the main components of the implementation, including the column generation, labeling algorithm, acceleration strategies, cutting planes, branching rules, and the three primal heuristics considered in the computational study.
The general principles of the pricing-for-integrality heuristic were presented in Section~\ref{sec:HeuristicPricing}, here, we describe the specific implementation details and choices made for the computational experiments.

\subsection{Column Generation}
To solve the ESPPRC that occurs as the pricing problem for the EVRPTW-CTS a bidirectional labeling algorithm \citep{righini2006symmetry} with the dynamic halfway-point strategy of \citet{tilk2017asymmetry} is implemented.
The resource extension functions and dominance criteria builds on the work of \cite{desaulniers2016exact} and their labeling algorithm for the ESSPRC associated with the EVRPTW, but with adaptions for handling the time slots at chargers.
A full description of resource extension functions and dominance criteria is omitted here for brevity, but can be found in Appendix \ref{sec:DetailedLabelingAlgorithm}.

Several strategies are used to accelerate the labeling algorithm.
First, the ng-route relaxation of \cite{baldacci2011new} relaxes elementarity by allowing certain node revisits while preventing revisits to nodes contained in predefined ng-sets.
This can substantially accelerate the labeling algorithm while still producing mostly elementary routes, although its performance depends on the size and composition of the ng-sets.
In our implementation, the ng-set of each node contains the node itself and its 15 nearest neighbours, including both customer and charging time-slot nodes.
Each charging time-slot node represents a unique charging opportunity and can therefore be visited at most once, unlike charging stations in most other EVRP variants.

Second, resource-based bounds are used to discard labels that cannot be extended into feasible routes.
A label is discarded if the earliest possible arrival at the depot exceeds the upper bound of the depot time window or if the remaining battery capacity is insufficient to reach either a charging station or the depot.

Third, heuristic pricing is performed on a subgraph of the original graph.
For each non-depot node \(i\in V\cup F\), the outgoing arcs are ordered by reduced cost, and only the five cheapest arcs to customers and the two cheapest arcs to charging nodes are retained.
The selected arcs are updated in each column-generation iteration based on current reduced costs, following \cite{desaulniers2008tabu}.

The acceleration strategies are applied sequentially.
Pricing is first performed heuristically using all accelerations and a reduced-cost threshold of \(-1\).
This threshold ensures that the heuristic stages prioritize columns with the potential to provide a meaningful improvement to the RMP objective.
If no column below the threshold is found, the graph reduction is relaxed to retain the 10 cheapest customer arcs and the three cheapest charging arcs from each node.
If this is still unsuccessful, graph reduction is disabled, while the ng-route relaxation and resource-based bounds remain active.
At this stage, any elementary column with negative reduced cost is accepted.
If the best negative-reduced-cost route is non-elementary, the ng-set size is increased by five and pricing is repeated, up to two times and a maximum ng-set size of 25.
If the route remains non-elementary, the ng-route relaxation is disabled and pricing is performed exactly with resource-based bounds remaining active.

In addition to the acceleration strategies described above, the dual values used in the pricing problem are stabilized.
We apply the smoothing method of \cite{neame2000nonsmooth}, where the stabilized dual values are computed as \(\widetilde{\pi}=\alpha\pi+(1-\alpha)\pi_{\mathrm{prev}}\).
Here, \(\pi\) denotes the current RMP dual solution, \(\pi_{\mathrm{prev}}\) denotes the stabilized dual values from the preceding iteration, and \(\alpha\in[0,1]\) controls the degree of stabilization.
We use $\alpha=0.5$ and apply stabilization in the root node until the average absolute change in the dual values between two consecutive iterations is below \(0.1\).

Before solving the pricing problem, the graph is preprocessed to remove infeasible arcs.
An arc is removed if its energy consumption exceeds the battery capacity.
An arc is also removed if the earliest possible arrival at its destination, when departing from the origin as early as possible, exceeds the upper bound of the destination time window.

\subsection{Cuts}
To strengthen the linear relaxation, we implement the 2-path cuts introduced for the VRPTW by \cite{kohl19992}.
For a subset of customers that cannot be served by a single vehicle, these cuts require at least two vehicles to be used.
The cuts are separated heuristically using a procedure similar to that of \cite{desaulniers2016exact}.
Given a fractional RMP solution, a support graph is constructed containing all nodes and the arcs with positive flow.
An exhaustive search is done to find all subsets of connected nodes in this graph, for which the total incoming flow is less than two.
To limit the computational effort, only subsets containing at most 10 nodes are considered.

For each candidate subset, we determine whether at least two vehicles are required.
If its total demand exceeds the vehicle capacity, this follows immediately.
Otherwise, an ESPPRC is solved on a graph containing the subset, the charging nodes, and the depot to determine whether all customers in the subset can be served by a single vehicle while satisfying the capacity, time-window, and battery constraints.
The ESPPRC is solved using the same labeling algorithm as the pricing problem, but with modified arc costs. Arcs entering customers in the candidate subset are assigned a cost of \(-1\), while all other arcs have zero cost.
Consequently, a route serving all customers in the subset has objective value equal to the negative number of customers in the subset.
If no such route exists, the corresponding 2-path cut is violated and added to the RMP.
In our implementation, 2-path cuts are separated only in the root node after the linear relaxation has been solved to optimality.
Preliminary experiments indicated that separating them elsewhere in the branch-and-bound tree generally did not pay off.

\subsection{Branching}
When the RMP solution is fractional, branching is performed hierarchically using two policies commonly applied in branch-and-price algorithms \citep[see, e.g.,][]{desaulniers2016exact,LAM2022105870,nafstad2025branch}.
First, we consider the total number of vehicles used. If this value is fractional, two child nodes are created: one where the number of vehicles is restricted to the nearest lower integer or less, and one where it is restricted to the nearest higher integer or more.
If the number of vehicles is integer, we instead consider the total flow on each arc.
If an arc \((i,j)\in A\) has fractional total flow, two child nodes are created: one where the arc is excluded and one where it is included.
For the child where the arc is excluded, arc \((i,j)\) is removed from the pricing graph.
For the child where the arc is included, all other outgoing arcs from node \(i\in V\cup F\) are removed.
The arc to branch on is simply selected as the first arc found with fractional total flow in the solution.
The branch-and-bound nodes are explored using the default SCIP settings.

\subsection{Primal Heuristics}
We consider three primal heuristics: a restricted master heuristic, a diving heuristic, and the proposed pricing for integrality heuristic.
The first two are executed after the LP relaxation has been solved and all cuts have been separated, while three root-node execution strategies are considered for the pricing for integrality heuristic, as described in Section~\ref{subsec:ImplementationDetailsPFHI}.
Beyond the root node, all three heuristics are applied every 100th processed node up to depth 20, provided that the global gap exceeds $0.5\%$ and the difference between the incumbent and the optimal RMP solution at the node is at least $50\%$ of the global gap.

\subsubsection{Restricted Master Heuristic}
The restricted master heuristic solves an integer problem over all generated columns using SCIP, with a time limit of 20 minutes.
Beyond the root node, the time limit is reduced to two minutes.

\subsubsection{Diving Heuristic}
The implemented diving heuristic is the diving heuristic with Limited Discrepancy Search (LDS) proposed by \cite{sadykov2019primal}.
Starting from the current solution of the RMP, the heuristic successively fixes fractional variables and re-solves the resulting RMP, thereby performing a dive toward an integer solution.
To introduce diversification and allow the heuristic to reconsider earlier fixing decisions, it relies on the limited discrepancy search paradigm introduced by \cite{harvey1995limited}.
LDS explores solutions that differ from the sequence of fixing decisions preferred by the diving rule, while limiting the number of such discrepancies.
This allows alternative branches of the heuristic search tree to be explored without performing an unrestricted search.
For fixing decisions, the diving heuristic uses a simple rule that selects the variable with the highest fractional value in the current solution of the RMP.

Following \cite{sadykov2019primal}, the maximum number of discrepancies is set to 3 and backtracking is allowed up to depth 2 in the diving search tree.
To speed up the heuristic, each residual master problem is solved with a time limit of 60 seconds, and only the heuristic pricing steps are used.
A related strategy is used by \cite{parmentier2023electric}, who perform pricing on a sparsified graph within their diving heuristic.
If heuristic pricing is unable to find a column with negative reduced cost, the pricing problem is not solved exactly and the current dive is terminated.
To improve solution quality, the remaining permitted LDS nodes are explored after the first feasible solution is found.
Beyond the root node, the heuristic performs a pure dive without backtracking, corresponding to a discrepancy limit of 0.

\subsubsection{Implementation of Pricing for Integrality Heuristic}\label{subsec:ImplementationDetailsPFHI}
The implementation of the proposed heuristic involves several choices, including when it should be applied during the BPC process.
As noted by \cite{maher2023integer}, its performance can depend considerably on the dual solution used, making the point of execution important.
To investigate this, we consider three strategies for applying the heuristic in the root node: (1) during pricing, (2) at an LP-optimal solution before separation, and (3) after separation of 2-path cuts from an LP-optimal solution.

Whenever the heuristic is applied, it constructs initial integer solutions from the current column pool, each serving as a separate starting point for destroy-and-repair.
For strategies (2) and (3), the heuristic is applied once and constructs five such solutions ($I_{\mathrm{start}}=5$).
For strategy (1), the heuristic is first applied when the dual values are considered to have stabilized, using the same criterion as for activating dual stabilization: a relative change in the dual values of less than (0.1\%).
It is then applied every fifth pricing iteration until LP optimality is reached or a maximum of five applications has been performed, constructing one initial solution ($I_{\mathrm{start}}=1$) each time.
Thus, strategy (1) constructs at most five initial solutions, while strategies (2) and (3) construct five each.
For each initial solution, the heuristic performs five destroy-and-repair iterations ($I_{\mathrm{destroy}}=5$).

For strategy (3), the restricted master heuristic is applied directly after the destroy-and-repair phase, as specified in the general design of the heuristic presented earlier.
For strategies (1) and (2), it is postponed until LP optimality and cut separation at the root node, allowing subsequent pricing to add columns.
In all three strategies, the final integer solution is thus constructed from the column pool available at the end of the root node, making the comparison fairer.
The restricted master component uses the same settings as the standalone restricted master heuristic, with a time limit of 20 minutes in the root node and two minutes elsewhere in the BPC tree.

When the heuristic is applied beyond the root node, it is run with $I_{\mathrm{start}}=1$ and $I_{\mathrm{destroy}}=5$ with the restricted master heuristic applied directly after the destroy-and-repair phase.
All active branching decisions are enforced both when constructing the initial solution and when solving the modified pricing problems during the repair phase.

In our implementation, the repair pricing problems are solved using all previously described acceleration strategies except for the ng-route relaxation.
Each repair pricing iteration is limited to 10 seconds and returns at most 10 columns.
If the time limit is reached, the pricing procedure is terminated and the columns found up to that point are used.
The ng-route relaxation is excluded because preliminary experiments showed that it reduced the performance of the heuristic.
In particular, it resulted in a substantial number of non-elementary columns, which appeared to be less useful for constructing high-quality integer solutions.

The settings specified above are fixed throughout the computational experiments.
The remaining parameters are \(\gamma\) and \(\beta\), which balance the influence of the dual values and heuristic weights in the modified pricing objective, and the destroy fraction \(\rho\), which determines the proportion of routes removed to create a partial solution for the repair process.
We believe these parameters are likely to have the most significant impact on the performance of the heuristic, as they directly control the guidance and extent of the destroy-and-repair procedure.
As such we conduct separate parameter analyses to examine their influence, with the results presented in Section~\ref{sec:Results}.

\section{Computational Results}\label{sec:Results}
In this section, we present computational results to evaluate the performance of the proposed Pricing-for-Integrality Heuristic (PFIH) as well as the impact of the charging time slots in the EVRPTW-CTS.
The experiments are conducted on adaptations of commonly used EVRPTW benchmark instances in which charging time slots have been introduced.
To facilitate further research, the generated instances and the code used in the experiments are publicly available at \url{https://gitlab.liu.se/lukev81/bpc-for-evrptw-cts-scip}.

\subsection{Instance Generation}
The instances are based on the EVRPTW dataset introduced by \cite{Schneider2014}, which in turn is derived from the widely used Solomon VRPTW benchmark instances \citep{Solomon1987}.
The Solomon dataset is divided into six instance classes: R1, R2, C1, C2, RC1, and RC2, which differ in terms of customer locations, time windows, and service requirements.
Out of these, we selected the C1, R1 and RC1 instances for our experiments, since they typically require more vehicles and therefore make the scheduling of charging more interesting.
Each EVRPTW instance contains 100 customers and 21 charging stations.

To obtain EVRPTW-CTS instances, the availability horizon of each physical charging station is divided into time slots of equal length.
The duration of each slot corresponds to the time required to charge \(80\%\) of the vehicle battery capacity.
Each available combination of a physical charging station and a time slot is represented by a distinct charging time-slot node in the network as described in Section~\ref{sec:ProblemFormulation}.
Each such node has capacity one, meaning that at most one vehicle can use the corresponding charging station during that time slot.
For each of the 29 original instances, we generate two variants with different levels of charging-slot availability: one in which approximately one-third of the charging time-slot nodes are available and one in which approximately half are available. This results in 58 instances, comprising 29 instances of each availability type.
Further details on the instance-generation procedure are provided in Appendix~\ref{sec:InstanceGeneration}.

\subsection{Results}
The algorithm was implemented in C++ using version 9.2.1 of the SCIP Optimization Suite \citep{bolusani2024scip}, with Gurobi 11.0.0 as the LP solver.
The implementation of the labeling algorithm was based on the code of \cite{eveborn2026subgraph}.
The experiments were conducted on a computing cluster equipped with Intel Xeon Gold 6130 processors.
Each run used 16 cores, a memory limit of 48 GB, and a SCIP time limit of two hours.
The reported gaps and running times are averages over three independent runs.

As described in Section~\ref{subsec:ImplementationDetailsPFHI}, the behavior of the heuristic depends on the weights \(\gamma\) and \(\beta\), the destroy fraction \(\rho\), and the point at which the heuristic is applied in the root node.
Our complete parameter analysis is presented in Appendix~\ref{sec:ParameterAnalysis}.
The performance of the heuristic was relatively insensitive to the choice of \(\gamma\) and \(\beta\), although larger values of \(\gamma\) generally increased solution times.
Based on these results, we use \(\gamma=0.6\) and \(\beta=0.4\).
The destroy fraction had a more clear effect on performance.
Small values of \(\rho\) generally resulted in larger optimality gaps, whereas large values increased solution times.
We therefore select \(\rho=0.6\) as a compromise between solution quality and computational effort.
The analysis of the execution point showed that applying the heuristic after cut separation produced the best overall performance.
One possible explanation is that the cuts strengthen the LP relaxation, leading to a different dual solution that may provide more useful guidance to the heuristic.
In the remaining experiments, PFIH is therefore applied after cut separation using \(\gamma=0.6\), \(\beta=0.4\), and \(\rho=0.6\).

We first evaluate the performance of PFIH in the root node of the BPC algorithm.
It is compared with a restricted master heuristic (RMH) and a diving heuristic, both of which represent commonly used primal-heuristic approaches in branch-and-price algorithms.
Table~\ref{tab:summary_rm_diving_pfih_after} reports the average root-node optimality gap and the average time spent in each heuristic.
Complete results for the individual instances are provided in Appendix~\ref{sec:FullResults}.
For each run, the optimality gap is calculated using its best primal solution at the end of the root node and a dual bound common to all runs of the same instance: the best root-node dual bound obtained across those runs.
Using a common dual bound ensures that differences in the reported gaps reflect differences in primal-solution quality rather than variations in the root-node dual bound between runs.
PFIH achieves the lowest average gap at both levels of charging-slot availability.
Across all 58 instances, its average root-node gap is 3.59\%, compared with 4.97\% for the diving heuristic and 8.81\% for the RMH.

\begin{table}
\centering
\footnotesize
\makebox[\linewidth][c]{%
\begin{threeparttable}
\caption{Root node optimality gap and heuristic time\label{tab:summary_rm_diving_pfih_after}}
\begin{tabular}{@{}l@{\quad}cc@{\quad}cc@{\quad}cc@{}}
\toprule
Series & \multicolumn{2}{c}{RMH} & \multicolumn{2}{c}{Diving} & \multicolumn{2}{c}{PFIH} \\
\cmidrule(lr){2-3} \cmidrule(lr){4-5} \cmidrule(lr){6-7}
 & Gap (\%) & Heur. time (s) & Gap (\%) & Heur. time (s) & Gap (\%) & Heur. time (s) \\ \midrule
\textbf{1/3-open} &  &  &  &  &  &  \\
C100-33 (9) & 3.38 & 87.12 & 2.64 & 69.25 & 2.39 & 177.70 \\
RC100-33 (8) & 7.64 & 217.39 & 6.91 & 107.74 & 4.68 & 175.70 \\
R100-33 (12) & 13.21 & 681.27 & 4.71 & 197.07 & 3.32 & 589.55 \\
\textbf{Average (29)} & \textbf{8.62} & \textbf{368.91} & \textbf{4.67} & \textbf{132.76} & \textbf{3.41} & \textbf{347.57} \\ \midrule
\textbf{1/2-open} &  &  &  &  &  &  \\
C100-50 (9) & 6.20 & 58.43 & 5.29 & 112.02 & 3.75 & 185.44 \\
RC100-50 (8) & 7.36 & 201.71 & 6.56 & 261.49 & 4.17 & 102.41 \\
R100-50 (12) & 12.21 & 653.76 & 4.37 & 300.35 & 3.52 & 451.55 \\
\textbf{Average (29)} & \textbf{9.00} & \textbf{344.30} & \textbf{5.26} & \textbf{231.18} & \textbf{3.77} & \textbf{272.65} \\ \midrule
\textbf{Overall average (58)} & \textbf{8.81} & \textbf{356.61} & \textbf{4.97} & \textbf{181.97} & \textbf{3.59} & \textbf{310.11} \\ \bottomrule
\end{tabular}
\begin{tablenotes}[flushleft]
\scriptsize\setlength{\labelsep}{0pt}
\item[] Average optimality gap and time spent in heuristic for the restricted master heuristic (RMH), the diving heuristic, and the pricing for integrality heuristic (PFIH) in the root node, grouped by instance series.
Number of instances considered in parentheses.
\end{tablenotes}
\end{threeparttable}}
\end{table}

As seen in the table, PFIH requires more time than the diving heuristic but less time than the RMH.
A considerable part of the PFIH running time is spent in its restricted master heuristic component, which is applied to the expanded column pool after the destroy-and-repair phase.
Its running time can therefore be reduced by imposing a shorter time limit on this component.
The corresponding analysis, presented in Appendix~\ref{sec:ParameterAnalysis}, shows that the running time can be reduced considerably with only a small deterioration in the resulting optimality gap.
Using a time limit for the restricted master heuristic component that results in an average total PFIH running time comparable to that of the diving heuristic, PFIH obtains an average optimality gap of approximately \(3.8\%\), compared to \(4.97\%\) for the diving heuristic.

We next evaluate the effect of the three heuristics on the complete BPC algorithm.
Table~\ref{tab:summary_full_runs} reports the average optimality gap after two hours and the number of instances solved to optimality.
Complete results for the individual instances are provided in Appendix~\ref{sec:FullResults}.
The results show that the advantage of PFIH over the two other heuristics persists after two hours of BPC, with an average optimality gap of \(2.63\%\) compared to \(3.87\%\) for the diving heuristic and \(6.63\%\) for the RMH.
When it comes to the number of instances solved to optimality, there is no clear advantage for any of the three heuristics.

\begin{table}
\centering
\footnotesize
\makebox[\linewidth][c]{%
\begin{threeparttable}
\caption{Optimality gap and instances solved to optimality\label{tab:summary_full_runs}}
\begin{tabular}{@{}l@{\quad}cc@{\quad}cc@{\quad}cc@{}}
\toprule
Series & \multicolumn{2}{c}{RMH} & \multicolumn{2}{c}{Diving} & \multicolumn{2}{c}{PFIH} \\
\cmidrule(lr){2-3} \cmidrule(lr){4-5} \cmidrule(lr){6-7}
 & Gap (\%) & Optimal (n) & Gap (\%) & Optimal (n) & Gap (\%) & Optimal (n) \\ \midrule
\textbf{1/3-open} &  &  &  &  &  &  \\
C100-33 (9) & 2.09 & 5 & 2.00 & 5 & 1.79 & 5 \\
RC100-33 (8) & 5.30 & 0 & 5.83 & 0 & 3.99 & 0 \\
R100-33 (12) & 7.03 & 4 & 3.49 & 3 & 2.03 & 4 \\
\textbf{Average (29)} & \textbf{5.02} & \textbf{9} & \textbf{3.67} & \textbf{8} & \textbf{2.50} & \textbf{9} \\ \midrule
\textbf{1/2-open} &  &  &  &  &  &  \\
C100-50 (9) & 2.79 & 3 & 3.65 & 3 & 2.43 & 3 \\
RC100-50 (8) & 6.45 & 0 & 5.92 & 0 & 3.46 & 0 \\
R100-50 (12) & 13.50 & 2 & 3.17 & 2 & 2.58 & 2 \\
\textbf{Average (29)} & \textbf{8.23} & \textbf{5} & \textbf{4.08} & \textbf{5} & \textbf{2.77} & \textbf{5} \\ \midrule
\textbf{Overall average (58)} & \textbf{6.63} & \textbf{14} & \textbf{3.87} & \textbf{13} & \textbf{2.63} & \textbf{14} \\ \bottomrule
\end{tabular}
\begin{tablenotes}[flushleft]
\scriptsize\setlength{\labelsep}{0pt}
\item[] Average optimality gap and number of instances solved to optimality for the restricted master (RMH), the diving heuristic, and the pricing for integrality heuristic (PFIH) after 2 hours of computation, grouped by instance series.
Number of instances considered in parentheses.
\end{tablenotes}
\end{threeparttable}}
\end{table}

Finally, we examine how the introduced charging-slot restrictions affect the instance solutions.
For this purpose, we compare the costs of the best primal solutions found for the \(1/3\)-open and \(1/2\)-open variants with those of the corresponding original EVRPTW instances.
The costs of the best primal solutions found are, on average, \(10.10\%\) higher for the \(1/3\)-open instances and \(4.26\%\) higher for the \(1/2\)-open instances than for the corresponding original EVRPTW instances.
These results indicate that the charging-slot restrictions have a meaningful impact on the resulting routes, with a greater effect when fewer charging slots are available.
Detailed results are provided in Appendix~\ref{sec:GapComparison}.

\section{Conclusion}\label{sec:Conclusion}
The proposed pricing-for-integrality primal heuristic appears to succeed in enriching the column pool with routes likely to contribute to high-quality integer solutions.
It uses a destroy-and-repair procedure to generate new routes, then applies a restricted master heuristic to the expanded column pool.
In computational experiments on the EVRPTW-CTS instances, the proposed heuristic reduced the average root-node optimality gap by nearly 60\% relative to a restricted master heuristic.
This suggests that the routes generated during repair make the column pool more useful for constructing high-quality integer solutions.
In a comparison with the implemented diving heuristic at similar average running times, the proposed heuristic achieved an approximately 25\% lower root-node gap.
After two hours of BPC, the reductions in average optimality gap relative to both baselines remained of similar magnitude.

The EVRPTW-CTS instances show that limited charging-slot availability has a meaningful effect on solution costs, with average increases of roughly 4--10\% relative to the corresponding EVRPTW instances.
Future work could evaluate the proposed heuristic on other vehicle routing problems to assess whether the results observed for the EVRPTW-CTS extend to other settings.

\section*{Acknowledgments}
The computations were enabled by resources provided by the National Academic Infrastructure for Supercomputing in Sweden (NAISS), partially funded by the Swedish Research Council through grant agreement no. 2022-06725.
We would also like to thank Scania and the Condore project group for fruitful discussions.


\newpage
\section*{Supplementary Material}
\appendix
\numberwithin{table}{section}
\numberwithin{figure}{section}
\fancyhead[L]{\scriptsize Appendix to \textbf{Eveborn and R\"{o}nnberg}: \textit{A Pricing-for-Integrality Heuristic for B\&P: An Application to the EVRPTW-CTS}}
\fancyhead[R]{}
\setstretch{1.1}
\fontsize{10}{12}\selectfont
\captionsetup{font=small}
\titlespacing*{\section}{0pt}{2.5ex plus 1ex minus .2ex}{1.2ex plus .2ex}
\setlength{\intextsep}{10pt plus 2pt minus 2pt}
\setlength{\textfloatsep}{12pt plus 2pt minus 4pt}
\setlength{\floatsep}{10pt plus 2pt minus 2pt}

\section{Compact Model}\label{sec:CompactModel}
The necessary notation for the Electric Vehicle Routing Problem with Time Windows and Charging Time Slots (EVRPTW-CTS) compact formulation is summarized in Table \ref{tab:notation}.
\begin{table}[ht]
\centering
\caption{Notation and model components}
\label{tab:notation}
\small
\scriptsize
\makebox[\linewidth][c]{%
\begin{tabular}{llll}
\toprule
Symbol & Description & Symbol & Description \\
\midrule
$0, n+1$ & Depot nodes.
  & $Q$ & Vehicle capacity. \\

$V$ & Set of customers.
  & $B$ & Battery capacity. \\

$F$ & Set of chargers duplicated by time discretization.
  & $g$ & Recharging rate. \\

$N$ & Union of all nodes.
  & $h$ & Energy consumption rate. \\

$K$ & Set of vehicles.
  && \\

$d_{ij}$ & Distance between nodes $i \in N,j \in N$.
    && \underline{\textbf{Variables}} \\

$t_{ij}$ & Travel time between nodes $i \in N,j \in N$.
    & $y_i^k$ & Battery level upon arrival at node $i \in N$. \\

  $e_i$ & Earliest start of service or charging at node $i \in N$.
    & $Y_i^k$ & Battery level upon departure from node $i \in N$. \\

$l_i$ & Latest start of service or charging at node $i \in N$.
& $p_i^k$ & Arrival time of vehicle $k \in K$ at node $i \in N$. \\

$s_i$ & Service time at node $i \in N$.
  & $x_{ij}^k$ & Binary variable equal to 1 if vehicle $k \in K$ travels\\

  $q_i$ & Demand at node $i \in N$.
  &&the arc between $i \in N$ and $j \in N$, 0 otherwise. \\
\bottomrule
\end{tabular}}
\end{table}

The compact formulation for EVRPTW-CTS is then given as follows:
{\thickmuskip=2mu plus 2mu
\medmuskip=1mu plus 1mu minus 1mu
\begin{subequations}\label{eq:compact_model}
\begin{align}
\text{minimize} \quad
    & \sum_{i \in N \setminus \{n+1\}} \sum_{j \in N \setminus \{0\},\, j \neq i} \sum_{k \in K} d_{ij} x_{ij}^k \label{main_problem:objective}, \\
\text{subject to} \quad
    & \sum_{j \in N \setminus \{0\},\, j \neq i} \sum_{k \in K} x_{ij}^k \ge 1, && i \in V \label{main_problem:customer_visits}, \\
    & \sum_{j \in N \setminus \{0\},\, j \neq i} \sum_{k \in K} x_{ij}^k \le 1, && i \in F \label{main_problem:charging_stations_visit}, \\
    & \sum_{i \in N \setminus \{n+1\},\, i \neq j} x_{ij}^k = \sum_{i \in N \setminus \{0\},\, i \neq j} x_{ji}^k, && j \in N,\ k \in K \label{main_problem:node_balance}, \\
    & \sum_{j \in N \setminus \{0, n+1\}} x_{0j}^k \le 1, && k \in K \label{main_problem:leave_station_once}, \\
    & \sum_{i \in V} \sum_{j \in N \setminus \{0\},\, j \neq i} q_i x_{ij}^k \le Q, && k \in K \label{main_problem:vehicle_capacity}, \\
    & y_j^k \le y_i^k - (h \cdot d_{ij}) x_{ij}^k + B(1-x_{ij}^k), && i \in V,\ j \in N \setminus \{0\},\, j \neq i,\ k \in K \label{main_problem:battery_on_arriving_from_customer}, \\
    & y_j^k \le Y_i^k - (h \cdot d_{ij}) x_{ij}^k + B(1-x_{ij}^k), && i \in F \cup \{0\},\ j \in N \setminus \{0\},\, j \neq i,\ k \in K \label{main_problem:battery_on_arriving_from_charger}, \\
    & y_i^k \le Y_i^k \le B, && i \in F \cup \{0\} \label{main_problem:max_battery}, \\
    & p_i^k + (t_{ij} + s_i)x_{ij}^k \le p_j^k + l_0(1-x_{ij}^k), && i \in V \cup \{0\},\ j \in N \setminus \{0\},\, j \neq i \label{main_problem:time_feasibility_leaving_from_customer}, \\
    & p_i^k + t_{ij} x_{ij}^k + g(Y_i^k-y_i^k) \le p_j^k + (l_0 + g B)(1-x_{ij}^k), && i \in F,\ j \in N \setminus \{0\},\, j \neq i,\ k \in K \label{main_problem:time_feasibility_leaving_from_charger}, \\
    & e_i \le p_i^k \le l_i, && i \in N \label{main_problem:time_windows_customers}, \\
    & p_i^k + g(Y_i^k-y_i^k) \le l_i, && i \in F,\ k \in K \label{main_problem:time_windows_chargers}, \\
    & x_{ij}^k \in \{0,1\}, && i \in N \setminus \{n+1\},\ j \in N \setminus \{0\},\, j \neq i,\ k \in K, \\
    & y_i^k \ge 0,\ Y_i^k \ge 0,\ p_i^k \ge 0, && i \in N,\ k \in K.
\end{align}
\end{subequations}}

The notation is inspired by the formulation in \cite{KUCUKOGLU2021107650}.
Some explanation of the objective function and constraints is provided below.
The objective function, given in \eqref{main_problem:objective}, minimizes the total distance traveled.
Constraint \eqref{main_problem:customer_visits} requires that each customer is visited at least once, while constraint \eqref{main_problem:charging_stations_visit} ensures that each charging time slot is visited at most once.
Node balance is enforced by \eqref{main_problem:node_balance}, and \eqref{main_problem:leave_station_once} ensures that each vehicle leaves the depot at most once.
Vehicle capacity is handled by \eqref{main_problem:vehicle_capacity}.
Battery consumption and recharging are modeled by constraints \eqref{main_problem:battery_on_arriving_from_customer}–\eqref{main_problem:max_battery}.
Time feasibility is ensured by constraints \eqref{main_problem:time_feasibility_leaving_from_customer}–\eqref{main_problem:time_windows_chargers}, where \eqref{main_problem:time_windows_chargers} enforces a hard time limit on the charging stations’ time windows.
Note that this differs from the time windows for customers, where it only restriction the start of service.
Compared to a standard EVRPTW formulation, constraint \eqref{main_problem:charging_stations_visit} is adapted and constraint \eqref{main_problem:time_windows_chargers} is added.

\section{Labeling Algorithm Details}\label{sec:DetailedLabelingAlgorithm}
In this appendix, the full labeling algorithm used for solving the pricing problem is described.
The resource extension functions, dominance criteria and notation described below follows the one presented in \cite{desaulniers2016exact} to a large extent, but with necessary modifications to handle the specifics of the EVRPTW-CTS.

The resources for a label $L_i$ of a partial path ${0} \rightarrow i \in N$ are defined as follows:
\begin{equation*}
\begin{aligned}
  T^{\text{cost}}_i &\colon
    \parbox[t]{0.8\textwidth}{Cost of the partial path.} \\
  T^{\text{load}}_i &\colon
    \parbox[t]{0.8\textwidth}{Load of the vehicle upon arrival at node $i$.} \\
  T^{\text{tMin}}_i       &\colon
    \parbox[t]{0.8\textwidth}{Earliest service start time at node $i$ assuming minimum recharge at visited chargers, while ensuring battery feasibility.} \\
  T^{\text{tMax}}_i       &\colon
    \parbox[t]{0.8\textwidth}{Earliest service start time at node $i$ assuming maximum recharge at visited chargers, while ensuring time-window feasibility.} \\
  T^{\text{rtMax}}_i      &\colon
    \parbox[t]{0.8\textwidth}{Maximum possible recharging time at node $i$ assuming minimum recharge at visited chargers, while ensuring battery feasibility.} \\
  (T^{\text{cust}_n}_i)_{n \in V} &\colon
    \parbox[t]{0.8\textwidth}{Indicates if the customer $n$ cannot be visited in the extension of the partial path.
    Is set to 1 if the customer cannot be visited, and 0 otherwise. A customer is considered unvisitable if either it is already visited or it cannot be reached because of battery, capacity or time window constraints.} \\
\end{aligned}
\end{equation*}

To define the resource extension functions, we first define some additional notation.
In addition to the parameters defined in appendix \ref{sec:CompactModel}, let $h_{ij} = d_{ij} \cdot h/g$ be the time needed to recharge the energy consumed when traveling from node $i$ to node $j$, and let $H = B/g$ be the time it takes to fully recharge the battery.
Also, let the travel time from node $i \in N$ to node $j \in N$ also include the service time at node $i$ moving forward, $t_{ij} = t_{ij} + s_i$.
The resource extension functions for extending a label $L_i$ at node $i \in N$ to a node $j \in N$ by traversing arc $(i,j) \in A$ are then defined as follows:
\begin{subequations}
  \begin{align}
    T_j^{\text{cost}} = & \; T_i^{\text{cost}} + c_{ij}, \\
    T_j^{\text{load}} = & \; T_i^{\text{load}} + q_j, \\
    T_j^{\text{tMin}} = & \; \text{max}(e_j, T_i^{\text{tMin}}+t_{ij}) + X_{ij}(T_i^{\text{tMin}}, T_i^{\text{tMax}}, T_i^{\text{rtMax}}), \\
    T_j^{\text{tMax}} = &
    \begin{cases}
        \text{min}(l_j, \text{max}(e_j, \text{min}(T_i^{\text{tMin}}+T_i^{\text{rtMax}}+t_{ij}, l_i + t_{ij}))), & \text{if } i \in F,\\
        \text{min}(l_j, \text{max}(e_j, T_i^{\text{tMax}}+t_{ij})), & \text{otherwise},
    \end{cases}\label{eq:Tj_tMax} \\
    T_j^{\text{rtMax}} = & \; \text{min}(H, \text{max}(0, T_i^{\text{rtMax}}-S_{ij}(T_i^{\text{tMin}}, T_i^{\text{tMax}}, T_i^{\text{rtMax}}) + h_{ij})), \\
    S_{ij} =  &
    \begin{cases}
        \text{max}(0, \text{min}(e_j-T_i^{\text{tMin}}-t_{ij}, \text{min}(T_i^{\text{rtMax}}, l_i - T_i^{\text{tMin}}))), & \text{if } i \in F,\\
        \text{max}(0, \text{min}(e_j - T_i^{\text{tMin}} - t_{ij}, T_i^{\text{tMax}} - T_i^{\text{tMin}})), & \text{otherwise},
    \end{cases}\label{eq:S_ij} \\
    T_j^{\text{cust}_n} = &
    \begin{cases}
        1 + T_i^{\text{cust}_n}, & \text{if } n = j,\\
         \text{max}(T_i^{\text{cust}_n},R(T_j^{\text{load}}, T_j^{\text{tMin}})), & \text{if } n \in N \setminus \{j\}.
    \end{cases}
  \end{align}
\end{subequations}
Where $X_{ij}(T_i^{\text{tMin}}, T_i^{\text{tMax}}, T_i^{\text{rtMax}}) = \text{max}(0, \text{max}(0, T_i^{\text{rtMax}}-S_{ij}(T_i^{\text{tMin}}, T_i^{\text{tMax}}, T_i^{\text{rtMax}})+h_{ij}-H))$, and $R(T_j^{\text{load}}, T_j^{\text{tMin}})$ is a function that returns 1 if $T_j^{\text{load}} + q_n > Q$ or $T_j^{\text{tMin}} + t_{jn} > l_n$, and 0 otherwise.

To account for the time slots at charger stations changes, in comparison to the uncapacitated case, changes are made to the calculation of the resource $T_j^{\text{tMax}}$ and the helper resource $S_{ij}$.
The helper resource $S_{ij}$ calculates the available slack time when extending a label from node $i$ to node $j$.
For $T_j^{\text{tMax}}$, the earliest start given maximum recharge, when extending from a charger station $i \in F$ to a node $j \in N$, it needs to consider if it is the maximum recharge, $T_i^{\text{tMin}}+T_i^{\text{rtMax}}+t_{ij}$, or the hard time window constraint, $l_i + t_{ij}$, that is limiting first.
This is done by taking the minimum of the two, which is seen in \eqref{eq:Tj_tMax}.
For the slack time $S_{ij}$, when extending from a charger station $i \in F$ to a node $j \in N$, it needs to be considered if it is the maximum recharge, $T_i^{\text{rtMax}}$, or the hard time window constraint, $l_i - T_i^{\text{tMin}}$, that limits the available slack time.
This is done by taking the minimum of the two, as seen in \eqref{eq:S_ij}.

An extension from node $i$ to node $j$ is feasible only if the following conditions are satisfied:
\begin{subequations}
  \begin{align}
    T_j^{\text{load}} \leq & Q, \\
    T_j^{\text{tMin}} \leq & l_j, \\
    T_j^{\text{tMin}} \leq & T_j^{\text{tMax}}, \\
    T_j^{\text{rtMax}} \leq & H, \\
    T_j^{\text{cust}_n} \leq & 1, ~\forall n \in N.
  \end{align}
\end{subequations}

For a label $L'_i$ at node $i$ to be dominated by another label $L_i$ at the same node, the following conditions must hold:
\begin{subequations}
  \begin{align}
    T_i^{\text{cost}} \leq & T_i^{'\text{cost}}, \\
    T_i^{\text{load}} \leq & T_i^{'\text{load}}, \\
    T_i^{\text{tMin}} \leq & T_i^{'\text{tMin}}, \\
    T_i^{\text{cust}_n} \leq & T_i^{'\text{cust}_n}, \forall n \in N, \\
    T_i^{\text{rtMax}} - (T_i^{\text{tMax}} - T_i^{\text{tMin}}) \leq & T_i^{'\text{rtMax}} - (T_i^{'\text{tMax}} - T_i^{'\text{tMin}}), \\
    T_i^{\text{rtMax}} - (T_i^{'\text{tMin}} - T_i^{'\text{tMin}}) \leq & T_i^{'\text{rtMax}}.
  \end{align}
\end{subequations}

For the backward labeling algorithm for a label $L_i$ of a partial path $i \in N \rightarrow {n+1}$, the resources are defined in the same manner as for the forward labeling algorithm, but for the backward direction.
They are denoted as $W^{\text{cost}}_i$, $W^{\text{load}}_i$, $W^{\text{tMin}}_i$, $W^{\text{tMax}}_i$, $W^{\text{rtMax}}_i$ and $(W^{\text{cust}_n}_i)_{n \in V}$.
For the initial label at the depot, i.e. $n+1$, all resources are set to zero, except for $W^{\text{tMin}}_{n+1}$ and $W^{\text{tMax}}_{n+1}$ which are set to the latest possible return time at the depot.

The resource extension functions for extending label $W_j$ at node $j \in N$ to a node $i \in N$ by traversing edge $(i,j) \in A$ are then defined as follows:
\begin{subequations}
  \begin{align}
    W_i^{\text{cost}} = & \; W_j^{\text{cost}} + c_{ij}, \\
    W_i^{\text{load}} = & \; W_j^{\text{load}} + q_j, \\
    W_i^{\text{tMin}} = & \;\min(l_i, W_j^{\text{tMin}}-t_{ij}) - \hat{X}_{ij}(W_j^{\text{tMin}}, W_j^{\text{tMax}}, W_j^{\text{rtMax}}), \\
    W_i^{\text{tMax}} =
    &\begin{cases}
      \max(e_i,
      \min(l_i,
      \max(W_j^{\text{tMin}}-W_j^{\text{rtMax}}-t_{ij}, e_j-t_{ij}))), & \text{if } j \in F,\\
      \max(e_i, \min(l_i, W_j^{\text{tMax}}-t_{ij})), & \text{otherwise},
    \end{cases} \label{eq:Wj_tMax}\\
    W_i^{\text{rtMax}} = & \;
      \min(H,
      \max(0,
        W_j^{\text{rtMax}}-\hat{S}_{ij}(W_j^{\text{tMin}}, W_j^{\text{tMax}}, W_j^{\text{rtMax}})
        +h_{ij}
      )), \\
    \hat{S}_{ij} =
    &\begin{cases}
      \max(
        0,
      \min(
        W_j^{\text{tMin}}-t_{ij}-l_i,
      \min(W_j^{\text{rtMax}}, W_j^{\text{tMin}}-e_j)
      )
      ),
        & \text{if } j \in F,\\
      \max(
        0,
      \min(
        W_j^{\text{tMin}}-t_{ij}-l_i,
        W_j^{\text{tMin}}-W_j^{\text{tMax}}
      )
      ),
        & \text{otherwise},
    \end{cases} \\
    W_i^{\text{cust}_n} =
    &\begin{cases}
        1 + W_j^{\text{cust}_n}, & \text{if } n = j,\\
         \text{max}(W_j^{\text{cust}_n}, \hat{R}(W_i^{\text{load}}, W_i^{\text{tMin}})), & \text{if } n \in N \setminus \{i\}.
    \end{cases}
  \end{align}
\end{subequations}

Where $\hat{X}_{ij}(W_j^{\text{tMin}}, W_j^{\text{tMax}}, W_j^{\text{rtMax}}) = \max(0, W_j^{\text{rtMax}} - \hat{S}_{ij}(W_j^{\text{tMin}}, W_j^{\text{tMax}}, W_j^{\text{rtMax}}) + h_{ij} - H)$ and $\hat{R}(W_i^{\text{load}}, W_i^{\text{tMin}})$ is a function that returns 1 if $W_i^{\text{load}} + q_n + q_i > Q$ or $W_i^{\text{tMin}} - t_{ni} < e_n$, and 0 otherwise.

To account for the time slots at charger stations changes, in comparison to the uncapacitated case, changes are made to the calculation of the resource $W_i^{\text{tMax}}$ and the helper resource $\hat{S}_{ij}$.
For $W_i^{\text{tMax}}$, the latest start given maximum recharge, when extending from a charger station $j \in F$ to a node $i \in N$, it needs to consider if it is the maximum recharge, $W_j^{\text{tMin}}-W_j^{\text{rtMax}}-t_{ij}$, or the hard time window constraint, $e_j-t_{ij}$, that is limiting first.
This is done by taking the maximum of the two, which is seen in \eqref{eq:Wj_tMax}.
For the slack time $\hat{S}_{ij}$, when extending from a charger station $j \in F$ to a node $i \in N$, it needs to be considered if it is the maximum recharge, $W_j^{\text{rtMax}}$, or the hard time window constraint, $W_j^{\text{tMin}}-e_j$, that limits the available slack time.

The extension from node $j$ to node $i$ is feasible if the following conditions are satisfied:
\begin{subequations}
  \begin{align}
    W_i^{\text{load}} \leq &~ Q, \\
    W_i^{\text{tMin}} \geq &~ e_i, \\
    W_i^{\text{tMin}} \geq &~ W_i^{\text{tMax}}, \\
    W_i^{\text{rtMax}} \leq &~ H, \\
    W_i^{\text{cust}_n} \leq &~ 1, ~\forall n \in N.
  \end{align}
\end{subequations}

A label $W'_i$ at node $i$ is dominated by another label $W_i$ at the same node if the following conditions are satisfied:
\begin{subequations}
  \begin{align}
    W_i^{\text{cost}} \leq &~ W_i^{'\text{cost}}, \\
    W_i^{\text{load}} \leq &~ W_i^{'\text{load}}, \\
    W_i^{\text{tMin}} \geq &~ W_i^{'\text{tMin}}, \\
    W_i^{\text{cust}_n} \leq &~ W_i^{'\text{cust}_n}, \forall n \in N, \\
    W_i^{\text{rtMax}} - (W_i^{\text{tMin}} - W_i^{\text{tMax}}) \leq &~ W_i^{'\text{rtMax}} - (W_i^{'\text{tMin}} - W_i^{'\text{tMax}}), \\
    W_i^{\text{rtMax}} - (W_i^{\text{tMin}} - W_i^{'\text{tMin}}) \leq &~ W_i^{'\text{rtMax}}.
  \end{align}
\end{subequations}

Finally, for the merge conditions for the bidirectional search, a forward label $W_i$ at node $i$ and a backward label $W_i$ at the same node can be merged to form a complete path, with cost $T_i^{\text{cost}}+W_i^{\text{cost}}$ from the depot to the depot if the following conditions are satisfied:
\begin{subequations}
  \begin{align}
    T_i^{\text{load}} + W_i^{\text{load}} \leq &~ Q, \\
    T_i^{\text{cust}_n} + W_i^{\text{cust}_n} \leq &~ 1, ~\forall n \in N \setminus \{i\}, \\
    T_i^{\text{tMin}} + \text{max}(0, T_i^{\text{rtMax}}+W_i^{\text{rtMax}+H}) \leq &~ W_i^{\text{tMin}} ~\text{if} ~ i \in V, \\
    \text{max}(0, T_i^{\text{rtMax}}+W_i^{\text{rtMax}+H}) \leq &~ (T_i^{\text{tMax}} - T_i^{\text{tMin}}) + (W_i^{\text{tMin}} - W_i^{\text{tMax}}), ~\text{if} ~ i \in F.
  \end{align}
\end{subequations}

\section{Instance Generation Details}\label{sec:InstanceGeneration}

The instances are based on the EVRPTW instances of \cite{Schneider2014}.
For each charging station in the original instance, denoted by $F'$, multiple charging time slots are generated in two steps: first, all potential time slots are created, and then a subset is designated as open.

To avoid unusable slots, only time slots reachable within the instance time horizon are considered.
Since the earliest allowed departure time from the depot is $e_0=0$ in all instances, the time horizon is $[0,l_0]$, where $l_0$ is the latest return time to the depot.
For each charger $i\in F'$, the time available for charging slots is therefore $t_i^{\text{available}}=l_0-2t_{0i}$, where $t_{0i}$ is the travel time between the depot and charger $i$.
The slot length is set to $t^{\text{slot}}=0.8B/g$, corresponding to the time required to charge 80\% of the battery capacity, where $B$ is the battery capacity and $g$ is the charging rate.

The number of potential slots at charger $i\in F'$ is $nr_i^{\text{chg}}=\left\lfloor t_i^{\text{available}}/t^{\text{slot}}\right\rfloor$.
The starting time of the first slot is $t^{\text{start}}*{i,1}=t_{0i}+\operatorname{random}\left(0:t_i^{\text{available}}-nr_i^{\text{chg}}t^{\text{slot}}\right)$, where $\operatorname{random}(a:b)$ returns a uniformly distributed integer between $a$ and $b$.
This randomly distributes the unused slack across chargers.
The remaining slots are placed consecutively according to $t^{\text{start}}*{i,l}=t^{\text{start}}*{i,l-1}+t^{\text{slot}}$, $l=2,\ldots,nr_i^{\text{chg}}$.

Once all potential slots have been created, a subset is randomly selected to be open.
Two instance types are defined by setting the number of open slots at each charger to $nr_i^{\text{open}}=\left\lceil0.33nr_i^{\text{chg}}\right\rceil$ for type 1 and $nr_i^{\text{open}}=\left\lceil0.50nr_i^{\text{chg}}\right\rceil$ for type 2.
Type 2 is constructed by adding slots to the corresponding type 1 instance and is therefore strictly more relaxed.

\section{Parameter Analysis}\label{sec:ParameterAnalysis}
\begin{figure}[ht]
    \centering
    \includegraphics[width=0.88\textwidth]{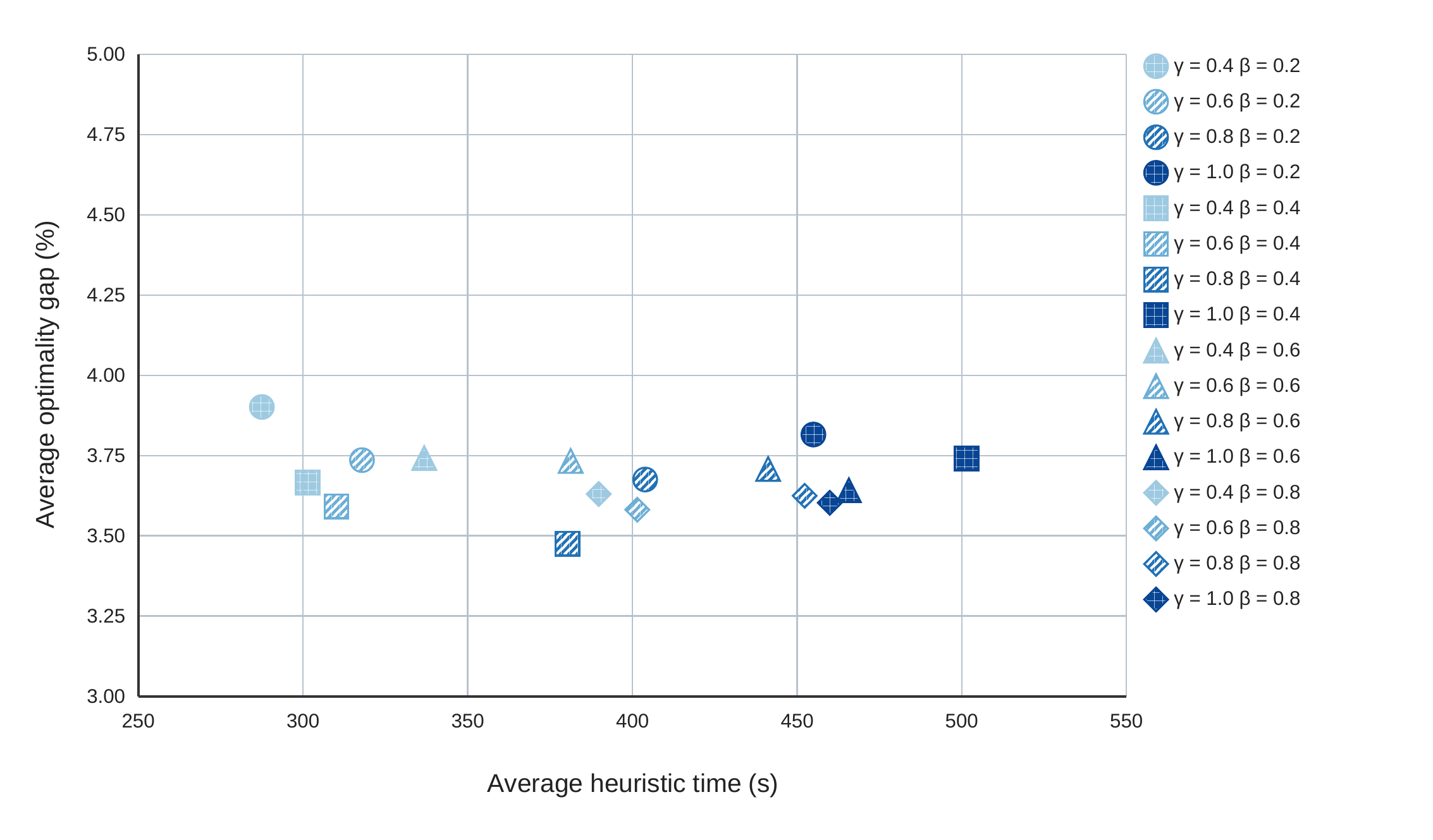}
    \caption{The figure shows the average optimality gap in the root node and the average time spent in the heuristic for different parameter settings of \(\gamma\) and \(\beta\).
    The heuristic was applied after cuts were separated, and the destroy ratio was set to \(\rho = 0.6\).}
    \label{fig:gamma_beta_parameter_analysis}
\end{figure}

\begin{figure}[ht]
    \centering
    \includegraphics[width=0.88\textwidth]{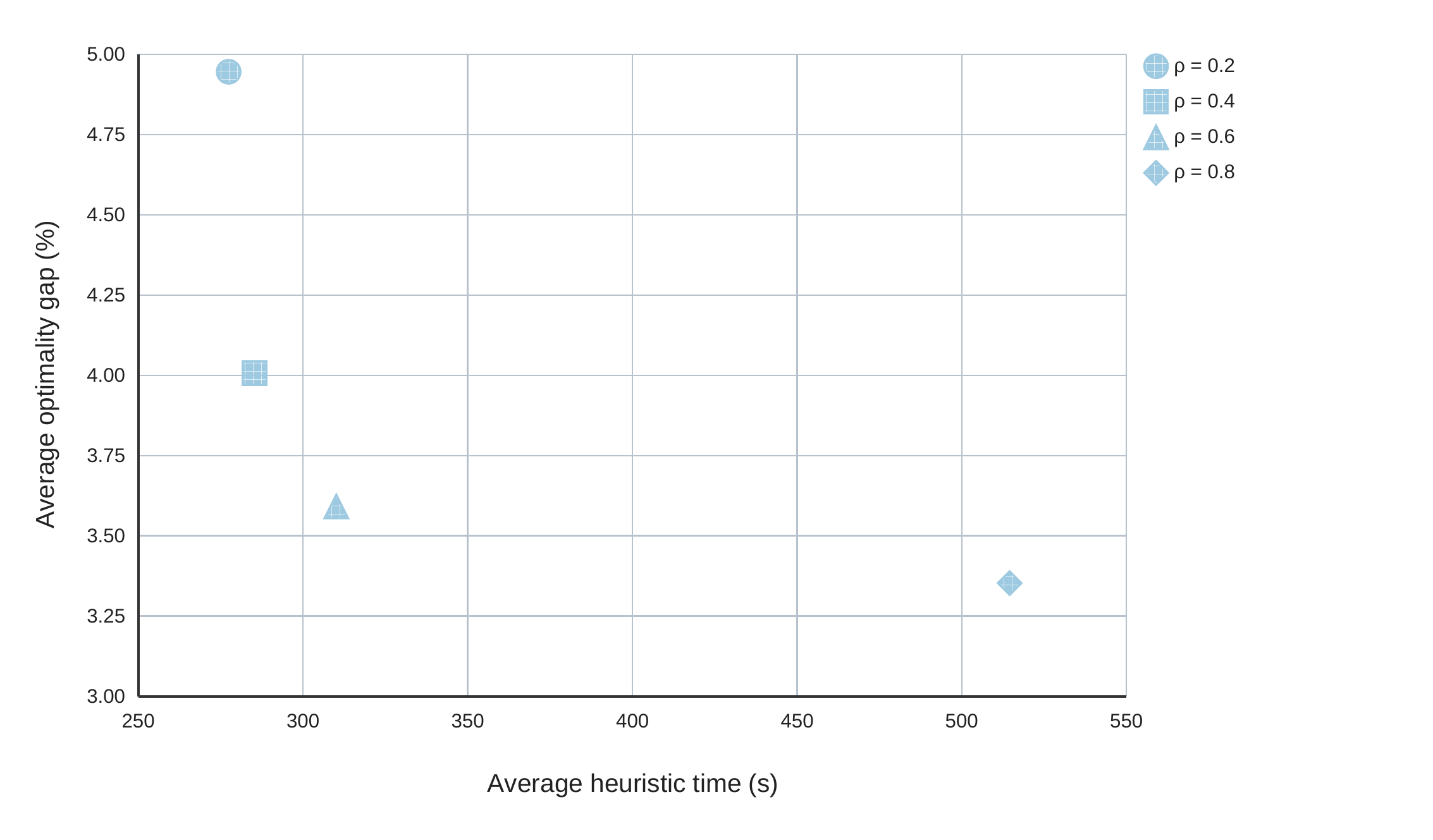}
    \caption{The figure shows the average optimality gap in the root node and the average time spent in the heuristic for different parameter settings of the destroy ratio \(\rho\).
    The heuristic was applied after cuts were separated, and the parameters \(\gamma\) and \(\beta\) were set to 0.6 and 0.4, respectively.}
    \label{fig:destroy_parameter_analysis}
\end{figure}

\begin{figure}[ht]
    \centering
    \includegraphics[width=0.88\textwidth]{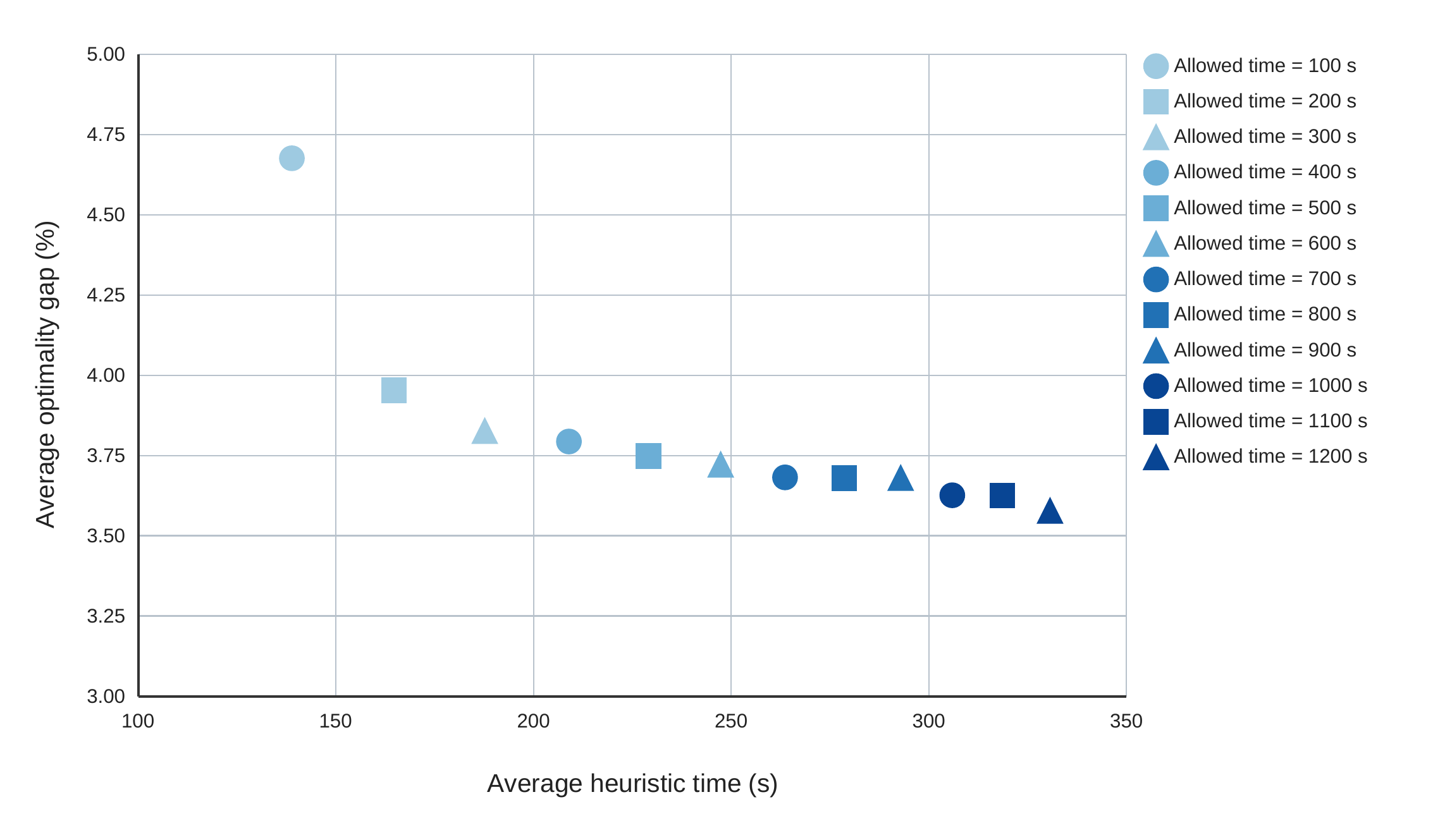}
    \caption{The figure shows the average optimality gap in the root node and the average time spent in the heuristic for different time limits of the restricted master heuristic part of PFIH.
    The heuristic was applied after cuts were separated, and the parameters \(\gamma\), \(\beta\), and \(\rho\) were set to 0.6, 0.4, and 0.6, respectively.}
    \label{fig:rm_time_limit_gap}
\end{figure}

This appendix presents the parameter analysis for the proposed pricing for integrality heuristic (PFIH).
Figures~\ref{fig:gamma_beta_parameter_analysis} and~\ref{fig:destroy_parameter_analysis} examine the weights \(\gamma\) and \(\beta\) and the destroy ratio \(\rho\), respectively.
Table~\ref{tab:summary_pfih} compares three points at which PFIH can be applied at the root node, and Figure~\ref{fig:rm_time_limit_gap} examines the time limit for the restricted master heuristic part of the heuristic.
All reported results are averages over three runs.

\begin{table}[ht]
\centering
\footnotesize
\makebox[\linewidth][c]{%
\begin{threeparttable}
\caption{Results for PFIH configurations\label{tab:summary_pfih}}
\scriptsize
\begin{tabular}{@{}l@{\quad}cc@{\quad}cc@{\quad}cc@{}}
\toprule
Series & \multicolumn{2}{c}{PFIH-during} & \multicolumn{2}{c}{PFIH-before} & \multicolumn{2}{c}{PFIH-after} \\
\cmidrule(lr){2-3} \cmidrule(lr){4-5} \cmidrule(lr){6-7}
 & Gap (\%) & Heur. time (s) & Gap (\%) & Heur. time (s) & Gap (\%) & Heur. time (s) \\ \midrule
\textbf{1/3-open} &  &  &  &  &  &  \\
C100-33 (9) & 2.65 & 124.39 & 2.67 & 99.88 & 2.39 & 177.70 \\
RC100-33 (8) & 5.16 & 188.98 & 4.76 & 103.10 & 4.68 & 175.70 \\
R100-33 (12) & 3.36 & 658.59 & 3.02 & 637.26 & 3.32 & 589.55 \\
\textbf{Average (29)} & \textbf{3.64} & \textbf{363.26} & \textbf{3.39} & \textbf{323.13} & \textbf{3.41} & \textbf{347.57} \\ \midrule
\textbf{1/2-open} &  &  &  &  &  &  \\
C100-50 (9) & 4.50 & 133.75 & 4.59 & 242.96 & 3.75 & 185.44 \\
RC100-50 (8) & 5.81 & 129.42 & 8.67 & 155.80 & 4.17 & 102.41 \\
R100-50 (12) & 3.73 & 645.85 & 3.05 & 589.01 & 3.52 & 451.55 \\
\textbf{Average (29)} & \textbf{4.54} & \textbf{344.46} & \textbf{5.08} & \textbf{362.11} & \textbf{3.77} & \textbf{272.65} \\ \midrule
\textbf{Overall average (58)} & \textbf{4.09} & \textbf{353.86} & \textbf{4.24} & \textbf{342.62} & \textbf{3.59} & \textbf{310.11} \\ \bottomrule
\end{tabular}
\begin{tablenotes}[flushleft]
\scriptsize\setlength{\labelsep}{0pt}
\item[] Comparison of three strategies for applying the pricing-for-integrality heuristic at the root node: (1) during pricing (PFIH-during), (2) before cut separation (PFIH-before), and (3) after cut separation (PFIH-after).
The table reports average optimality gaps and time spent in the heuristic, grouped by instance series; the number of instances in each series is given in parentheses.
The parameters \(\gamma\), \(\beta\), and \(\rho\) are set to 0.6, 0.4, and 0.6, respectively.
If no primal solution is found in a run, that run contributes a 100\% gap to the average.
\end{tablenotes}
\end{threeparttable}}
\end{table}

\clearpage

\section{Gap Comparison}\label{sec:GapComparison}
This appendix compares the costs of the best primal solutions found for the \(1/3\)-open and \(1/2\)-open variants with those for the corresponding original EVRPTW instances.
The percentage gap is calculated as \(100(z_{\mathrm{variant}}-z_{\mathrm{original}})/z_{\mathrm{original}}\), where \(z_{\mathrm{variant}}\) and \(z_{\mathrm{original}}\) denote the cost of the best primal solution found for the variant and the original instance, respectively.
The results are reported in Table~\ref{tab:gap_comparison}.
\begin{table}[ht]
\centering
\footnotesize
\makebox[\linewidth][c]{%
\begin{threeparttable}
\caption{Comparison of Gaps\label{tab:gap_comparison}}
\scriptsize
\setlength{\tabcolsep}{5pt}
\begin{tabular}{@{}l@{\hspace{14pt}}r@{\hspace{14pt}}rr@{\hspace{14pt}}rr@{}}
\toprule
 & \multicolumn{1}{c}{Standard} & \multicolumn{2}{c}{1/3-open} & \multicolumn{2}{c}{1/2-open} \\
\cmidrule(lr){3-4} \cmidrule(lr){5-6}
Instance & Primal & Primal & Gap (\%) & Primal & Gap (\%) \\ \midrule
c101 & 1043.44 & 1140.80 & 9.33 & 1104.59 & 5.86 \\
c102 & 1011.82 & 1091.57 & 7.88 & 1063.21 & 5.08 \\
c103 & 993.57 & 1073.73 & 8.07 & 1051.46 & 5.83 \\
c104 & 898.13 & 1052.72 & 17.21 & 1026.86 & 14.33 \\
c105 & 1015.86 & 1111.79 & 9.44 & 1078.86 & 6.20 \\
c106 & 1008.07 & 1089.26 & 8.05 & 1058.47 & 5.00 \\
c107 & 1012.25 & 1081.48 & 6.84 & 1061.22 & 4.84 \\
c108 & 1004.61 & 1049.50 & 4.47 & 1025.07 & 2.04 \\
c109 & 936.37 & 1040.78 & 11.15 & 1013.00 & 8.18 \\
r101 & 1600.93 & 1800.36 & 12.46 & 1634.77 & 2.11 \\
r102 & 1420.76 & 1591.30 & 12.00 & 1462.06 & 2.91 \\
r103 & 1206.52 & 1374.86 & 13.95 & 1277.03 & 5.84 \\
r104 & 1073.72 & 1130.11 & 5.25 & 1089.77 & 1.49 \\
r105 & 1347.77 & 1498.81 & 11.21 & 1388.06 & 2.99 \\
r106 & 1256.94 & 1395.72 & 11.04 & 1288.16 & 2.48 \\
r107 & 1113.56 & 1204.60 & 8.18 & 1132.29 & 1.68 \\
r108 & 1060.42 & 1092.84 & 3.06 & 1060.42 & 0.00 \\
r109 & 1172.95 & 1271.23 & 8.38 & 1206.38 & 2.85 \\
r110 & 1076.80 & 1122.20 & 4.22 & 1089.57 & 1.19 \\
r111 & 1075.82 & 1140.04 & 5.97 & 1094.07 & 1.70 \\
r112 & 1007.87 & 1078.77 & 7.03 & 1051.53 & 4.33 \\
rc101 & 1647.15 & 1940.61 & 17.82 & 1730.63 & 5.07 \\
rc102 & 1503.18 & 1709.02 & 13.69 & 1602.33 & 6.60 \\
rc103 & 1301.42 & 1505.62 & 15.69 & 1381.19 & 6.13 \\
rc104 & 1210.64 & 1312.26 & 8.39 & 1297.70 & 7.19 \\
rc105 & 1439.58 & 1708.35 & 18.67 & 1493.10 & 3.72 \\
rc106 & 1387.09 & 1577.85 & 13.75 & 1445.99 & 4.25 \\
rc107 & 1265.60 & 1396.99 & 10.38 & 1274.06 & 0.67 \\
rc108 & 1197.34 & 1310.20 & 9.43 & 1231.69 & 2.87 \\
\midrule
\textbf{Average} &  &  & \textbf{10.10} &  & \textbf{4.26} \\
\bottomrule
\end{tabular}
\end{threeparttable}}
\end{table}

\section{Full Results}\label{sec:FullResults}
The tables in this appendix report results from two sets of experiments for all instances: runs terminated after the root node and full branch-and-bound runs.
All input data were rounded to two decimal places.
The reported results are averages over three runs.
For each set of experiments, results are reported separately for the restricted master heuristic (RMH), the diving heuristic, and the proposed pricing-for-integrality heuristic (PFIH).
Although the time limit was set to two hours, SCIP occasionally exceeded it slightly while completing internal processes.
In these cases, the reported time is 7,200 seconds.
If no primal solution was found in one or two of the three runs, the primal bound for each such run was set to twice its dual bound before calculating the averages and gaps.
Primal bounds affected by this convention are marked with an asterisk.
Table~\ref{tab:detailed_results_root} reports the dual bound, primal bound, optimality gap, total time, and time spent in the heuristic for the root-node runs.
To make the optimality gaps comparable across settings in the root node, they are calculated using the best dual bound obtained across all settings and runs for each instance.
Table~\ref{tab:detailed_full_runs} reports the primal bound, dual bound, optimality gap, total time, and total number of nodes for the full branch-and-bound runs.

\begin{table}
\centering
\footnotesize
\makebox[\linewidth][c]{%
\begin{threeparttable}
\caption{Root node results for RMH, Diving, and PFIH\label{tab:detailed_results_root}}
\scriptsize
\setlength{\tabcolsep}{2pt}
\begin{tabular}{@{}lr@{\hspace{10pt}}rrrr@{\hspace{10pt}}rrrr@{\hspace{10pt}}rrrr@{}}
\toprule
 &  & \multicolumn{4}{c}{RMH} & \multicolumn{4}{c}{Diving} & \multicolumn{4}{c}{PFIH} \\
\cmidrule(lr){3-6} \cmidrule(lr){7-10} \cmidrule(lr){11-14}
Instance & Dual & Primal & Gap & Time & H-time & Primal & Gap & Time & H-time & Primal & Gap & Time & H-time \\
 &  &  & (\%) & (s) & (s) &  & (\%) & (s) & (s) &  & (\%) & (s) & (s) \\ \midrule
c101-33 & 1140.80 & 1140.80 & 0.00 & 1.94 & 0.00 & 1140.80 & 0.00 & 1.62 & 0.00 & 1140.80 & 0.00 & 1.87 & 0.00 \\
c102-33 & 1087.70 & 1094.57 & 0.63 & 19.06 & 1.09 & 1091.57 & 0.36 & 38.93 & 18.19 & 1092.72 & 0.46 & 29.90 & 11.40 \\
c103-33 & 1050.10 & 1107.64 & 5.48 & 152.20 & 15.23 & 1090.04 & 3.80 & 252.29 & 111.30 & 1091.52 & 3.94 & 180.72 & 50.88 \\
c104-33 & 963.16 & 1075.69 & 11.68 & 610.86 & 5.59 & 1066.27 & 10.71 & 1012.63 & 299.43 & 1059.18 & 9.97 & 1120.60 & 507.49 \\
c105-33 & 1111.79 & 1111.79 & 0.00 & 6.06 & 0.00 & 1111.79 & 0.00 & 6.23 & 0.00 & 1111.79 & 0.00 & 6.28 & 0.00 \\
c106-33 & 1081.33 & 1090.05 & 0.81 & 11.70 & 1.15 & 1089.60 & 0.77 & 35.68 & 25.69 & 1089.43 & 0.75 & 24.67 & 13.59 \\
c107-33 & 1067.97 & 1098.83 & 2.89 & 17.16 & 4.12 & 1087.28 & 1.81 & 56.32 & 43.18 & 1082.27 & 1.34 & 32.34 & 15.35 \\
c108-33 & 1043.70 & 1076.47 & 3.14 & 162.14 & 131.19 & 1062.28 & 1.78 & 93.77 & 61.13 & 1058.47 & 1.42 & 64.07 & 28.91 \\
c109-33 & 1010.59 & 1068.90 & 5.77 & 756.23 & 625.70 & 1056.80 & 4.57 & 189.87 & 64.35 & 1046.83 & 3.59 & 1102.02 & 971.70 \\
r101-33 & 1796.44 & 1811.15 & 0.82 & 10.00 & 2.62 & 1807.96 & 0.64 & 27.92 & 20.54 & 1802.07 & 0.31 & 18.73 & 11.25 \\
r102-33 & 1585.80 & 1608.14 & 1.41 & 24.16 & 8.05 & 1617.80 & 2.02 & 51.08 & 34.80 & 1595.70 & 0.62 & 32.74 & 12.71 \\
r103-33 & 1362.93 & 1408.53 & 3.35 & 490.00 & 431.71 & 1420.32 & 4.21 & 124.92 & 61.56 & 1391.88 & 2.12 & 397.19 & 337.60 \\
r104-33 & 1080.00 & 1257.04 & 16.39 & 1746.59 & 1200.11 & 1175.88 & 8.88 & 1224.15 & 474.49 & 1136.37 & 5.22 & 2218.82 & 1566.80 \\
r105-33 & 1487.62 & 1520.57 & 2.21 & 397.44 & 56.51 & 1516.85 & 1.96 & 387.52 & 74.37 & 1499.92 & 0.83 & 322.74 & 13.01 \\
r106-33 & 1378.21 & 1422.69 & 3.23 & 422.89 & 80.72 & 1427.22 & 3.56 & 354.05 & 39.84 & 1402.52 & 1.76 & 359.59 & 19.45 \\
r107-33 & 1167.99 & 1322.48 & 13.23 & 1573.88 & 1200.10 & 1242.94 & 6.42 & 484.87 & 100.94 & 1216.43 & 4.15 & 1627.53 & 1226.54 \\
r108-33 & 1046.08 & 1371.38 & 31.10 & 2141.76 & 1200.10 & 1120.51 & 7.12 & 1688.16 & 688.51 & 1106.81 & 5.81 & 2779.54 & 1757.74 \\
r109-33 & 1237.52 & 1327.21 & 7.25 & 1791.25 & 1200.07 & 1294.90 & 4.64 & 664.43 & 54.60 & 1278.54 & 3.32 & 1755.73 & 1215.17 \\
r110-33 & 1103.60 & 1191.66 & 7.98 & 1256.43 & 786.90 & 1178.77 & 6.81 & 760.49 & 174.77 & 1133.40 & 2.70 & 596.92 & 123.98 \\
r111-33 & 1122.10 & 1206.86 & 7.55 & 1410.65 & 808.31 & 1157.47 & 3.15 & 655.55 & 123.08 & 1143.55 & 1.91 & 763.37 & 88.43 \\
r112-33 & 1031.56 & *1691.63 & 63.99 & 1895.21 & 1200.08 & 1104.46 & 7.07 & 1641.84 & 517.31 & 1146.33 & 11.13 & 1297.86 & 701.90 \\
rc101-33 & 1878.31 & 1950.72 & 3.85 & 216.88 & 9.22 & 2000.03 & 6.48 & 209.71 & 24.60 & 1946.68 & 3.64 & 235.77 & 34.46 \\
rc102-33 & 1660.41 & 1755.67 & 5.74 & 242.84 & 9.14 & 1785.23 & 7.52 & 251.54 & 43.72 & 1720.59 & 3.62 & 328.37 & 52.12 \\
rc103-33 & 1455.11 & 1538.42 & 5.73 & 196.19 & 17.36 & 1530.21 & 5.16 & 283.11 & 82.83 & 1516.89 & 4.25 & 255.39 & 39.74 \\
rc104-33 & 1248.96 & 1424.10 & 14.02 & 1640.08 & 485.95 & 1343.60 & 7.58 & 1309.50 & 170.46 & 1326.85 & 6.24 & 1174.92 & 203.31 \\
rc105-33 & 1668.77 & 1755.63 & 5.21 & 441.54 & 28.71 & 1738.34 & 4.17 & 490.53 & 75.53 & 1721.64 & 3.17 & 376.20 & 39.57 \\
rc106-33 & 1543.30 & 1633.49 & 5.84 & 493.13 & 24.66 & 1616.71 & 4.76 & 563.65 & 41.37 & 1579.86 & 2.37 & 477.28 & 17.43 \\
rc107-33 & 1321.15 & 1473.30 & 11.52 & 1565.12 & 864.90 & 1465.19 & 10.90 & 919.76 & 162.60 & 1410.74 & 6.78 & 1180.68 & 442.02 \\
rc108-33 & 1230.71 & 1343.84 & 9.19 & 1182.61 & 299.19 & 1337.53 & 8.68 & 1204.72 & 260.79 & 1321.13 & 7.35 & 1448.68 & 576.98 \\
c101-50 & 1099.80 & 1123.24 & 2.13 & 8.49 & 1.05 & 1134.55 & 3.16 & 11.23 & 4.98 & 1120.71 & 1.90 & 12.75 & 6.34 \\
c102-50 & 1058.02 & 1072.18 & 1.34 & 85.17 & 9.09 & 1066.09 & 0.76 & 122.17 & 26.89 & 1064.23 & 0.59 & 78.90 & 17.00 \\
c103-50 & 1020.56 & 1096.11 & 7.40 & 343.17 & 58.20 & 1102.34 & 8.01 & 530.72 & 200.84 & 1057.63 & 3.63 & 388.76 & 115.84 \\
c104-50 & 940.18 & 1051.13 & 11.80 & 690.72 & 10.83 & 1036.74 & 10.27 & 1205.63 & 441.90 & 1022.34 & 8.74 & 1956.24 & 1292.25 \\
c105-50 & 1072.51 & 1124.13 & 4.81 & 18.47 & 3.09 & 1101.84 & 2.73 & 38.07 & 23.19 & 1109.03 & 3.41 & 30.42 & 16.16 \\
c106-50 & 1045.12 & 1097.75 & 5.04 & 25.90 & 2.34 & 1086.57 & 3.97 & 58.77 & 40.48 & 1089.79 & 4.27 & 34.64 & 15.31 \\
c107-50 & 1039.86 & 1099.73 & 5.76 & 33.12 & 4.91 & 1103.55 & 6.13 & 95.80 & 68.14 & 1079.97 & 3.86 & 47.02 & 19.33 \\
c108-50 & 1014.03 & 1067.22 & 5.25 & 96.82 & 7.66 & 1066.89 & 5.21 & 215.97 & 113.32 & 1032.63 & 1.83 & 124.91 & 36.66 \\
c109-50 & 973.40 & 1092.95 & 12.28 & 630.18 & 428.69 & 1045.32 & 7.39 & 256.60 & 88.48 & 1027.42 & 5.55 & 334.29 & 150.03 \\
r101-50 & 1626.75 & 1646.02 & 1.18 & 16.75 & 6.51 & 1642.74 & 0.98 & 31.62 & 21.96 & 1634.77 & 0.49 & 26.43 & 13.72 \\
r102-50 & 1449.11 & 1469.99 & 1.44 & 81.21 & 8.78 & 1468.73 & 1.35 & 151.03 & 71.19 & 1462.90 & 0.95 & 88.71 & 13.99 \\
r103-50 & 1256.14 & 1304.01 & 3.81 & 173.48 & 75.19 & 1319.03 & 5.01 & 202.38 & 120.51 & 1278.63 & 1.79 & 171.00 & 78.07 \\
r104-50 & 1049.65 & 1250.63 & 19.15 & 2187.73 & 1200.10 & 1109.06 & 5.66 & 1963.11 & 505.17 & 1092.79 & 4.11 & 2418.70 & 1616.22 \\
r105-50 & 1374.21 & 1402.90 & 2.09 & 271.01 & 130.43 & 1414.95 & 2.96 & 230.92 & 97.38 & 1388.89 & 1.07 & 121.50 & 20.71 \\
r106-50 & 1274.44 & 1334.66 & 4.73 & 1475.74 & 1200.07 & 1311.81 & 2.93 & 305.49 & 104.78 & 1294.85 & 1.60 & 267.40 & 45.66 \\
r107-50 & 1113.63 & 1157.50 & 3.94 & 640.52 & 353.64 & 1136.41 & 2.05 & 611.05 & 295.39 & 1139.87 & 2.36 & 741.96 & 495.23 \\
r108-50 & 1013.72 & 1339.54 & 32.14 & 2110.69 & 1200.10 & 1089.10 & 7.44 & 1866.84 & 807.96 & 1060.54 & 4.62 & 2207.21 & 1103.11 \\
r109-50 & 1183.00 & 1281.20 & 8.30 & 1910.50 & 1200.08 & 1238.69 & 4.71 & 719.38 & 99.56 & 1212.78 & 2.52 & 930.54 & 349.96 \\
r110-50 & 1061.45 & 1201.81 & 13.22 & 1366.13 & 1052.68 & 1120.32 & 5.55 & 640.13 & 295.63 & 1097.01 & 3.35 & 425.97 & 95.96 \\
r111-50 & 1073.06 & 1131.46 & 5.44 & 530.24 & 217.45 & 1117.17 & 4.11 & 662.08 & 314.10 & 1096.98 & 2.23 & 458.10 & 63.44 \\
r112-50 & 991.01 & *1496.81 & 51.04 & 2599.21 & 1200.09 & 1087.25 & 9.71 & 2956.40 & 870.58 & 1161.56 & 17.21 & 2987.91 & 1522.55 \\
rc101-50 & 1702.01 & 1755.82 & 3.16 & 260.93 & 8.56 & 1766.50 & 3.79 & 320.76 & 63.84 & 1732.87 & 1.81 & 266.45 & 21.34 \\
rc102-50 & 1564.99 & 1642.58 & 4.96 & 580.78 & 21.61 & 1672.00 & 6.84 & 700.76 & 88.49 & 1606.68 & 2.66 & 644.71 & 23.03 \\
rc103-50 & 1350.30 & 1466.74 & 8.62 & 414.89 & 20.90 & 1432.78 & 6.11 & 587.74 & 194.45 & 1396.53 & 3.42 & 420.71 & 45.21 \\
rc104-50 & 1188.75 & 1349.71 & 13.54 & 2553.87 & 464.89 & 1301.80 & 9.51 & 2875.48 & 512.00 & 1293.76 & 8.83 & 2795.31 & 301.38 \\
rc105-50 & 1464.59 & 1541.31 & 5.24 & 331.16 & 10.45 & 1538.10 & 5.02 & 471.10 & 124.81 & 1504.04 & 2.69 & 445.83 & 116.28 \\
rc106-50 & 1419.40 & 1491.66 & 5.09 & 891.30 & 216.35 & 1488.03 & 4.84 & 765.37 & 103.32 & 1459.91 & 2.85 & 881.87 & 131.41 \\
rc107-50 & 1227.50 & 1351.10 & 10.07 & 2211.23 & 741.11 & 1321.72 & 7.68 & 1991.59 & 279.16 & 1299.25 & 5.84 & 1645.45 & 64.18 \\
rc108-50 & 1172.29 & 1267.99 & 8.16 & 1512.42 & 129.82 & 1273.93 & 8.67 & 2248.39 & 725.83 & 1233.37 & 5.21 & 1701.25 & 116.47 \\
\bottomrule
\end{tabular}
\end{threeparttable}}
\end{table}

\begin{table}
\centering
\footnotesize
\makebox[\linewidth][c]{%
\begin{threeparttable}
\caption{Full-run results for RMH, Diving, and PFIH\label{tab:detailed_full_runs}}
\scriptsize
\setlength{\tabcolsep}{2pt}
\begin{tabular}{@{}l@{\hspace{8pt}}rrrrr@{\hspace{8pt}}rrrrr@{\hspace{8pt}}rrrrr@{}}
\toprule
 & \multicolumn{5}{c}{RMH} & \multicolumn{5}{c}{Diving} & \multicolumn{5}{c}{PFIH} \\
\cmidrule(lr){2-6} \cmidrule(lr){7-11} \cmidrule(lr){12-16}
Instance & Primal & Dual & Gap & Time & Nodes & Primal & Dual & Gap & Time & Nodes & Primal & Dual & Gap & Time & Nodes \\
 &  &  & (\%) & (s) &  &  &  & (\%) & (s) &  &  &  & (\%) & (s) &  \\ \midrule
c101-33 & 1140.80 & 1140.80 & 0.00 & 2.14 & 1 & 1140.80 & 1140.80 & 0.00 & 2.20 & 1 & 1140.80 & 1140.80 & 0.00 & 2.01 & 1 \\
c102-33 & 1091.57 & 1091.57 & 0.00 & 106.92 & 23 & 1091.57 & 1091.57 & 0.00 & 92.17 & 24 & 1091.57 & 1091.57 & 0.00 & 102.49 & 22 \\
c103-33 & 1092.01 & 1053.08 & 3.70 & 7200.00 & 644 & 1087.83 & 1052.86 & 3.32 & 7200.00 & 529 & 1086.91 & 1053.25 & 3.20 & 7200.00 & 701 \\
c104-33 & 1076.26 & 964.97 & 11.53 & 7200.00 & 93 & 1065.05 & 965.01 & 10.37 & 7200.00 & 143 & 1060.87 & 964.99 & 9.94 & 7200.00 & 183 \\
c105-33 & 1111.79 & 1111.79 & 0.00 & 7.93 & 1 & 1111.79 & 1111.79 & 0.00 & 6.23 & 1 & 1111.79 & 1111.79 & 0.00 & 7.69 & 1 \\
c106-33 & 1089.26 & 1089.26 & 0.00 & 373.96 & 595 & 1089.26 & 1089.26 & 0.00 & 364.22 & 618 & 1089.26 & 1089.26 & 0.00 & 319.93 & 579 \\
c107-33 & 1081.48 & 1077.23 & 0.39 & 7200.00 & 11327 & 1083.61 & 1076.77 & 0.64 & 7200.00 & 10174 & 1081.48 & 1077.10 & 0.41 & 7200.00 & 9764 \\
c108-33 & 1049.50 & 1049.50 & 0.00 & 382.88 & 208 & 1049.50 & 1049.50 & 0.00 & 270.19 & 259 & 1049.50 & 1049.50 & 0.00 & 690.68 & 183 \\
c109-33 & 1047.16 & 1015.07 & 3.16 & 7200.00 & 1834 & 1052.75 & 1015.27 & 3.69 & 7200.00 & 2397 & 1041.25 & 1015.15 & 2.57 & 7200.00 & 1673 \\
r101-33 & 1800.36 & 1800.36 & 0.00 & 18.63 & 53 & 1800.36 & 1800.36 & 0.00 & 31.00 & 53 & 1800.36 & 1800.36 & 0.00 & 28.67 & 52 \\
r102-33 & 1591.30 & 1591.30 & 0.00 & 198.96 & 227 & 1591.30 & 1591.30 & 0.00 & 187.34 & 188 & 1591.30 & 1591.30 & 0.00 & 211.24 & 207 \\
r103-33 & 1374.86 & 1374.86 & 0.00 & 6279.67 & 2595 & 1374.86 & 1373.05 & 0.13 & 6872.24 & 3106 & 1374.86 & 1374.86 & 0.00 & 5683.39 & 2313 \\
r104-33 & 1206.32 & 1083.88 & 11.30 & 7200.00 & 12 & 1162.58 & 1084.34 & 7.22 & 7200.00 & 13 & 1142.64 & 1084.33 & 5.38 & 7200.00 & 10 \\
r105-33 & 1498.81 & 1498.81 & 0.00 & 448.37 & 479 & 1498.81 & 1498.81 & 0.00 & 551.67 & 419 & 1498.81 & 1498.81 & 0.00 & 492.46 & 461 \\
r106-33 & 1395.72 & 1389.97 & 0.41 & 7200.00 & 3382 & 1419.89 & 1386.19 & 2.43 & 7200.00 & 2455 & 1395.72 & 1390.13 & 0.40 & 7200.00 & 3613 \\
r107-33 & 1297.64 & 1173.47 & 10.58 & 7200.00 & 75 & 1219.51 & 1173.48 & 3.92 & 7200.00 & 61 & 1207.46 & 1173.52 & 2.89 & 7200.00 & 66 \\
r108-33 & 1237.94 & 1045.92 & 18.36 & 7200.00 & 7 & 1132.78 & 1045.92 & 8.30 & 7200.00 & 7 & 1102.10 & 1045.92 & 5.37 & 7200.00 & 5 \\
r109-33 & 1301.18 & 1244.78 & 4.53 & 7200.00 & 663 & 1303.04 & 1244.45 & 4.71 & 7200.00 & 976 & 1275.90 & 1243.37 & 2.62 & 7200.00 & 737 \\
r110-33 & 1134.77 & 1112.25 & 2.03 & 7200.00 & 159 & 1175.65 & 1112.05 & 5.72 & 7200.00 & 146 & 1124.52 & 1112.45 & 1.09 & 7200.00 & 154 \\
r111-33 & 1212.30 & 1122.10 & 8.04 & 7200.00 & 7 & 1156.19 & 1122.10 & 3.04 & 7200.00 & 3 & 1142.15 & 1122.10 & 1.79 & 7200.00 & 3 \\
r112-33 & 1332.48 & 1031.56 & 29.17 & 7200.00 & 7 & 1097.37 & 1031.56 & 6.38 & 7200.00 & 7 & 1081.05 & 1031.56 & 4.80 & 7200.00 & 3 \\
rc101-33 & 1940.61 & 1892.84 & 2.52 & 7200.00 & 6759 & 1940.95 & 1897.32 & 2.30 & 7200.00 & 29770 & 1940.61 & 1892.07 & 2.57 & 7200.00 & 5683 \\
rc102-33 & 1721.80 & 1673.85 & 2.86 & 7200.00 & 3849 & 1774.25 & 1672.07 & 6.11 & 7200.00 & 5637 & 1711.05 & 1673.56 & 2.24 & 7200.00 & 2969 \\
rc103-33 & 1519.30 & 1462.34 & 3.90 & 7200.00 & 417 & 1527.02 & 1463.33 & 4.35 & 7200.00 & 536 & 1507.26 & 1463.88 & 2.96 & 7200.00 & 515 \\
rc104-33 & 1382.04 & 1251.78 & 10.41 & 7200.00 & 54 & 1341.85 & 1251.78 & 7.20 & 7200.00 & 61 & 1346.93 & 1251.87 & 7.59 & 7200.00 & 33 \\
rc105-33 & 1711.31 & 1688.91 & 1.33 & 7200.00 & 2758 & 1757.93 & 1687.83 & 4.15 & 7200.00 & 2775 & 1708.35 & 1690.59 & 1.05 & 7200.00 & 3425 \\
rc106-33 & 1577.85 & 1554.22 & 1.52 & 7200.00 & 1972 & 1592.02 & 1554.34 & 2.42 & 7200.00 & 2452 & 1577.85 & 1555.53 & 1.43 & 7200.00 & 2310 \\
rc107-33 & 1446.96 & 1323.67 & 9.31 & 7200.00 & 132 & 1472.23 & 1323.49 & 11.24 & 7200.00 & 145 & 1405.59 & 1322.64 & 6.27 & 7200.00 & 102 \\
rc108-33 & 1360.18 & 1230.71 & 10.52 & 7200.00 & 7 & 1339.90 & 1230.71 & 8.87 & 7200.00 & 7 & 1326.82 & 1230.71 & 7.81 & 7200.00 & 7 \\
c101-50 & 1104.59 & 1104.59 & 0.00 & 16.06 & 11 & 1104.59 & 1104.59 & 0.00 & 28.44 & 12 & 1104.59 & 1104.59 & 0.00 & 22.20 & 12 \\
c102-50 & 1063.21 & 1063.21 & 0.00 & 526.31 & 60 & 1063.21 & 1063.21 & 0.00 & 480.66 & 53 & 1063.21 & 1063.21 & 0.00 & 442.89 & 52 \\
c103-50 & 1067.07 & 1024.24 & 4.18 & 7200.00 & 294 & 1105.48 & 1024.41 & 7.91 & 7200.00 & 259 & 1052.96 & 1025.20 & 2.71 & 7200.00 & 397 \\
c104-50 & 1051.69 & 942.11 & 11.63 & 7200.00 & 96 & 1030.47 & 942.24 & 9.36 & 7200.00 & 75 & 1044.50 & 942.11 & 10.87 & 7200.00 & 54 \\
c105-50 & 1078.86 & 1078.86 & 0.00 & 659.60 & 1203 & 1078.86 & 1078.86 & 0.00 & 526.75 & 1225 & 1078.86 & 1078.86 & 0.00 & 434.47 & 1203 \\
c106-50 & 1062.48 & 1048.63 & 1.32 & 7200.00 & 3080 & 1075.34 & 1048.71 & 2.54 & 7200.00 & 3429 & 1061.07 & 1048.99 & 1.15 & 7200.00 & 3723 \\
c107-50 & 1063.47 & 1045.84 & 1.69 & 7200.00 & 5408 & 1067.07 & 1046.14 & 2.00 & 7200.00 & 9059 & 1071.38 & 1045.81 & 2.44 & 7200.00 & 4611 \\
c108-50 & 1025.73 & 1021.06 & 0.46 & 7200.00 & 3032 & 1058.16 & 1019.76 & 3.77 & 7200.00 & 3416 & 1025.64 & 1021.10 & 0.44 & 7200.00 & 2984 \\
c109-50 & 1035.51 & 978.43 & 5.83 & 7200.00 & 615 & 1048.97 & 978.05 & 7.25 & 7200.00 & 680 & 1019.72 & 978.40 & 4.22 & 7200.00 & 773 \\
r101-50 & 1634.77 & 1634.77 & 0.00 & 83.32 & 327 & 1634.77 & 1634.77 & 0.00 & 106.49 & 303 & 1634.77 & 1634.77 & 0.00 & 79.12 & 292 \\
r102-50 & 1462.06 & 1459.51 & 0.17 & 7200.00 & 5636 & 1462.73 & 1458.11 & 0.32 & 7200.00 & 5374 & 1462.06 & 1460.96 & 0.08 & 7103.79 & 5868 \\
r103-50 & 1277.03 & 1263.20 & 1.10 & 7200.00 & 1616 & 1277.03 & 1263.53 & 1.07 & 7200.00 & 1932 & 1277.03 & 1262.93 & 1.12 & 7200.00 & 1634 \\
r104-50 & 1239.20 & 1049.60 & 18.06 & 7200.00 & 11 & 1107.23 & 1049.60 & 5.49 & 7200.00 & 12 & 1114.33 & 1049.60 & 6.17 & 7200.00 & 9 \\
r105-50 & 1388.06 & 1388.06 & 0.00 & 1833.65 & 2775 & 1388.06 & 1388.06 & 0.00 & 1630.45 & 3230 & 1388.06 & 1388.06 & 0.00 & 1585.68 & 2641 \\
r106-50 & 1288.16 & 1284.52 & 0.28 & 7200.00 & 1483 & 1290.25 & 1283.91 & 0.49 & 7200.00 & 1534 & 1288.16 & 1284.87 & 0.26 & 7200.00 & 1530 \\
r107-50 & 1140.87 & 1120.17 & 1.85 & 7200.00 & 222 & 1135.83 & 1120.63 & 1.36 & 7200.00 & 249 & 1136.64 & 1120.51 & 1.44 & 7200.00 & 210 \\
r108-50 & 1202.76 & 1013.72 & 18.65 & 7200.00 & 7 & 1078.76 & 1013.77 & 6.41 & 7200.00 & 6 & 1079.35 & 1013.72 & 6.47 & 7200.00 & 7 \\
r109-50 & 1206.59 & 1190.17 & 1.38 & 7200.00 & 944 & 1242.64 & 1188.73 & 4.53 & 7200.00 & 1118 & 1206.66 & 1190.96 & 1.32 & 7200.00 & 974 \\
r110-50 & 1218.09 & 1065.89 & 14.28 & 7200.00 & 69 & 1106.09 & 1067.18 & 3.65 & 7200.00 & 76 & 1099.87 & 1067.61 & 3.02 & 7200.00 & 77 \\
r111-50 & 1146.23 & 1078.77 & 6.25 & 7200.00 & 50 & 1114.53 & 1080.30 & 3.17 & 7200.00 & 66 & 1098.20 & 1080.51 & 1.64 & 7200.00 & 78 \\
r112-50 & - & 991.01 & - & 7200.00 & 7 & 1105.19 & 991.01 & 11.52 & 7200.00 & 7 & 1084.61 & 991.01 & 9.44 & 7200.00 & 6 \\
rc101-50 & 1730.63 & 1710.03 & 1.20 & 7200.00 & 5988 & 1731.60 & 1713.17 & 1.08 & 7200.00 & 17760 & 1730.63 & 1709.23 & 1.25 & 7200.00 & 5343 \\
rc102-50 & 1613.16 & 1572.62 & 2.58 & 7200.00 & 1693 & 1674.72 & 1571.37 & 6.58 & 7200.00 & 1921 & 1602.49 & 1572.30 & 1.92 & 7200.00 & 1413 \\
rc103-50 & 1442.00 & 1356.41 & 6.31 & 7200.00 & 148 & 1432.18 & 1356.08 & 5.61 & 7200.00 & 296 & 1382.55 & 1356.29 & 1.94 & 7200.00 & 184 \\
rc104-50 & 1364.16 & 1188.75 & 14.76 & 7200.00 & 11 & 1323.96 & 1188.75 & 11.37 & 7200.00 & 11 & 1312.73 & 1188.75 & 10.43 & 7200.00 & 11 \\
rc105-50 & 1503.66 & 1480.57 & 1.56 & 7200.00 & 2138 & 1506.99 & 1479.83 & 1.84 & 7200.00 & 2952 & 1493.10 & 1480.91 & 0.82 & 7200.00 & 2317 \\
rc106-50 & 1450.59 & 1425.57 & 1.75 & 7200.00 & 1250 & 1495.87 & 1425.07 & 4.97 & 7200.00 & 1326 & 1449.14 & 1425.18 & 1.68 & 7200.00 & 1272 \\
rc107-50 & 1357.59 & 1230.93 & 10.29 & 7200.00 & 56 & 1324.43 & 1230.86 & 7.60 & 7200.00 & 84 & 1283.18 & 1230.03 & 4.32 & 7200.00 & 84 \\
rc108-50 & 1326.82 & 1172.29 & 13.18 & 7200.00 & 7 & 1269.68 & 1172.29 & 8.31 & 7200.00 & 7 & 1234.47 & 1172.29 & 5.30 & 7200.00 & 6 \\
\bottomrule
\end{tabular}
\end{threeparttable}}
\end{table}
\clearpage

\end{document}